\documentclass[a4paper,reqno,twoside]{amsart} 
\usepackage{amssymb}
\usepackage{latexsym}
\usepackage{amsmath}
\usepackage{mathrsfs}
\usepackage[all,2cell]{xy}
        \UseAllTwocells \SilentMatrices
\usepackage{ifthen}
\usepackage{comment}

\theoremstyle{plain}
\newtheorem{propn}{Proposition}[section]
\newtheorem{thm}[propn]{Theorem}
\newtheorem{lemma}[propn]{Lemma}
\newtheorem{cor}[propn]{Corollary}

\theoremstyle{definition}
\newtheorem*{defn}{Definition}

\theoremstyle{remark}
\newtheorem*{rem}{Remark}
\newtheorem*{rems}{Remarks}
\newtheorem*{notn}{Notation}

\newcommand{\ve}{\varepsilon}

\newcommand{\varphitilde}{\wt{\varphi}}
\newcommand{\varphiol}{\ol{\varphi}}
\newcommand{\phitilde}{\wt{\phi}}
\newcommand{\chitilde}{\wt{\chi}}
\newcommand{\psitilde}{\wt{\psi}}
\newcommand{\psiol}{\ol{\psi}}
\newcommand{\phiol}{\ol{\phi}}
\newcommand{\chiol}{\ol{\chi}}
\newcommand{\gammatilde}{\wt{\gamma}}

\newcommand{\Rbar}{\mathcal{R}_\sigma}
\newcommand{\Ebar}{\mathcal{E}_\sigma}
\newcommand{\Rmap}{\mathcal{R}}
\newcommand{\Emap}{\mathcal{E}}

\newcommand{\qg}{\mathbb{G}}

\newcommand{\Mil}{\mathsf{M}}
\newcommand{\Nil}{\mathsf{N}}
\newcommand{\Altilde}{\wt{\mathsf{A}}}
\newcommand{\Alol}{\ol{\mathsf{A}}}
\newcommand{\Al}{\mathsf{A}}
\newcommand{\Biltilde}{\wt{\mathsf{B}}}
\newcommand{\Bilol}{\ol{\mathsf{B}}}
\newcommand{\Bil}{\mathsf{B}}     
\newcommand{\Cil}{\mathsf{C}}
\newcommand{\Ciltilde}{\wt{\mathsf{C}}}

\newcommand{\Hil}{\mathsf{H}}
\newcommand{\hil}{\mathsf{h}}
\newcommand{\Kil}{\mathsf{K}}
\newcommand{\kil}{\mathsf{k}}
\newcommand{\Vil}{\mathsf{V}}

\newcommand{\Wil}{\mathsf{W}}

\newcommand{\Blg}{\mathcal{B}}
\newcommand{\blg}{\mathsf{B}}

\newcommand{\Exps}{\mathcal{E}}
\newcommand{\Fock}{\mathcal{F}}
\newcommand{\FFock}{\mathcal{F}}

\newcommand{\Adjointable}{\mathcal{L}}
\newcommand{\KAdjointable}{\mathcal{K}}

\newcommand{\Cstar}{C^*}

\newcommand{\coproduct}{\Delta}
\newcommand{\coproducttilde}{\wt{\Delta}}
\newcommand{\coproductol}{\ol{\Delta}}
\newcommand{\counit}{\epsilon}
\newcommand{\Cou}{\counit}
\newcommand{\counittilde}{\wt{\epsilon}}
\newcommand{\counitol}{\ol{\epsilon}}

\newcommand{\kilhat}{\wh{\kil}}
\newcommand{\hilhat}{\wh{\hil}}

\newcommand{\khat}{\wh{\kil}}

\newcommand{\chat}{\wh{c}}
\newcommand{\dhat}{\wh{d}}
\newcommand{\ghat}{\wh{g}}
\newcommand{\fhat}{\wh{f}}

\newcommand{\Proc}{\mathbb{P}}
\newcommand{\Procstar}{\mathbb{P}_\conv}

\newcommand{\Step}{\mathbb{S}}
\newcommand{\Real}{\mathbb{R}}
\newcommand{\Rplus}{\Real_+}
\newcommand{\Comp}{\mathbb{C}}

\newcommand{\bn}{\mathbb{N}}

\newcommand{\ZZ}{\mathbb{Z}}

\newcommand{\dQS}{\Delta^{\mathrm{QS}}}

\newcommand{\cb}{{\text{\tu{cb}}}}

\newcommand{\Bstrict}{B_{\beta}}
\newcommand{\Bsigma}{B_{\sigma}}
\newcommand{\Buw}{B_{\sigma}}
\newcommand{\CBstrict}{CB_{\beta}}
\newcommand{\CPstrict}{CP_{\beta}}
\newcommand{\CBuw}{CB_{\sigma}}
\newcommand{\CPuw}{CP_{\sigma}}

\newcommand{\Linol}{\ol{\Lin}}

\newcommand{\wh}{\widehat}
\newcommand{\wt}{\widetilde}
\newcommand{\ol}{\overline}
\newcommand{\ul}{\underline}
\newcommand{\ot}{\otimes}
\newcommand{\otM}{\otimes_{\mathrm{M}}}
\newcommand{\ottilde}{\,\wt{\otimes}\,}
\newcommand{\otol}{\,\ol{\otimes}\,}
\newcommand{\otul}{\,\ul{\otimes}\,}

\newcommand{\conv}{\star}
\newcommand{\convtilde}{\,\wt{\star}\,}
\newcommand{\convol}{\,\ol{\star}\,}

\newcommand{\la}{\langle}
\newcommand{\ra}{\rangle}

\newcommand{\fullcomp}{\boldsymbol{\cdot}}

\newcommand{\tu}{\textup}

\DeclareMathOperator{\Ran}{Ran}
\DeclareMathOperator{\Lin}{Lin}
\DeclareMathOperator{\Ker}{Ker}
\DeclareMathOperator{\id}{id}

\newenvironment{alist}
{

\begin{enumerate}}
{\end{enumerate}}

\newenvironment{rlist}
{

\begin{enumerate}}
{\end{enumerate}}

\numberwithin{equation}{section}
\begin{document}

\title
[Quantum stochastic convolution cocycles] {Quantum stochastic
convolution cocycles III}

\author[Lindsay]{J.\ Martin Lindsay}
\address{Department of Mathematics and Statistics,
Lancaster University, Lancaster LA1 4YF, U.K.}
\email{j.m.lindsay@lancaster.ac.uk}
\author[Skalski]{Adam G.\ Skalski}
\address{Department of Mathematics and Statistics,
Lancaster University, Lancaster LA1 4YF, U.K.}
\curraddr{Institute of Mathematics, Polish Academy of Sciences, ul.
\'Sniadeckich 8, 00-956 Warsaw, Poland.}

\email{a.skalski@impan.pl}


 \begin{abstract}
  Every Markov-regular quantum L\'evy process on a
 multiplier $C^*$-bialgebra is shown to be equivalent to one governed
 by a quantum stochastic differential equation, and the generating
 functionals of norm-continuous convolution semigroups on a
 multiplier $C^*$-bialgebra are then completely characterised.
 These results are achieved by extending the theory of quantum
 L\'{e}vy processes on a compact quantum group,
 and more generally quantum stochastic convolution cocycles on a
 $C^*$-bialgebra, to locally compact quantum groups and
 multiplier $C^*$-bialgebras. Strict extension results obtained by
 Kustermans, together with automatic strictness properties developed here, are
 exploited to obtain existence and uniqueness for coalgebraic quantum
 stochastic differential equations in this setting. Then, working in the
 universal enveloping von Neumann bialgebra, we characterise the stochastic
 generators of Markov-regular, *-homomorphic (respectively completely positive
 and contractive), quantum stochastic convolution cocycles.
 \end{abstract}


\maketitle

\section*{Introduction}
\label{Section: Introduction}

Let $\qg$ be a locally compact quantum semigroup. The main results
of this paper are as follows: a
 concrete realisation of each abstract \emph{quantum L\'evy
 process} on $\qg$ which is \emph{Markov-regular} (that is, has
 norm-continuous expectation semigroup)
 as a quantum stochastic process in some Fock space (Corollary 6.2),
 and a characterisation of the generators of norm-continuous
  convolution semigroups
  of states on $\qg$ (Theorem 6.3).
 These are achieved by the development of a general
 theory of quantum stochastic evolutions with
 tensor-independent identically distributed increments, in the
 broader context of \emph{locally compact quantum semigroups} (in
 other words multiplier $C^*$-bialgebras), culminating in
 infinitesimal characterisations of Markov-regular *-homomorphic,
 respectively completely positive and contractive, quantum
 stochastic convolution cocycles (Theorem~\ref{Theorem: cocycle}).
 As a consequence of our results, a large family of examples of
 quantum L\'evy processes on $\qg$ is generated. These are indexed
 by a nondegenerate representation of $\qg$ (as a C$^*$-algebra) and
 a vector from the representing Hilbert space.

The notion of quantum L\'evy process generalises that of classical
L\'evy process on a semigroup. It was first introduced by Accardi,
Sch\"urmann and von Waldenfels, in the purely algebraic framework of
$^*$-bialgebras (\cite{asw}), and was further developed by
Sch\"urmann and others (\cite{Schurmann}, \cite{Uwe}) who, in
particular, extended it to other noncommutative forms of
independence (free, boolean and monotone), still in the algebraic
context. Inspired by Sch\"urmann's reconstruction theorem, which
states that every quantum L\'evy process on a $^*$-bialgebra can be
equivalently realised on a symmetric Fock space, we first showed how
the algebraic theory of quantum L\'evy processes can be extended to
the natural setting of \emph{quantum stochastic convolution
cocycles} (\cite{QSCC1}). These are families of linear maps
$(l_t)_{t\geq 0}$ from a $^*$-bialgebra $\Blg$ to operators on the
symmetric Fock space $\FFock$, over a Hilbert space of the form
$L^2(\Rplus; \kil)$, satisfying the following cocycle identity with
respect to the ampliated CCR flow $(\sigma_t)_{t\geq 0}$:
\[
l_{s+t} = l_s \star (\sigma_s \circ l_t), \;\;\; s,t \geq 0,
\]
together with regularity and adaptedness conditions. Our approach
enabled us to then establish a theory of quantum L\'evy processes on
\emph{compact quantum groups} and, more generally, quantum
stochastic convolution cocycles on operator space coalgebras
(\cite{QSCC2}).
 An extensive list of examples is included in Section 8 of that
 paper.
 The recent development of a satisfactory theory of \emph{locally
compact} quantum groups (\cite{kuv}) provides the challenge which is
addressed in the current work, namely to extend our analysis to the
locally compact realm.

 On the algebraic level the
theories of quantum stochastic convolution cocycles on compact and
locally compact quantum semigroups look similar, however their
analytic aspects have a rather different nature. Whereas the
coproduct on a compact quantum semigroup $\blg$ takes values in the
spatial tensor product $\blg \ot \blg$, which led us to an
operator-space theoretic development of the theory and enabled us to
establish the main results in the corresponding natural category of
operator-space coalgebras, a \emph{noncompact} locally compact
quantum semigroup $\blg$ is a nonunital $C^*$-algebra whose
coproduct takes values in the multiplier algebra of $\blg \ot \blg$.
Consequently, $C^*$-algebraic methods come more to the fore, with
the strict topology (\cite{Lance}), strict maps (\cite{Kus}) and
enveloping von Neumann algebras all playing crucial roles. The
modern approach to quantum stochastics
involves matrix spaces as a natural tool for combining
$C^*$-algebraic quantum state spaces with von Neumann algebraic
quantum noise (\cite{LWexistence}), and this was successfully
exploited in~\cite{QSCC2}. In the context of the present paper, the
strict topology on the initial $C^*$-algebra has to be harnessed to
the matrix-space technology, so that both may be exploited in
tandem.

The paper divides into three parts. The first part of the paper
(Sections ~\ref{Section: Extensions} and~\ref{Section: Bialgebras})
is devoted to a careful analysis of relations between the strict
topology on multiplier algebras and the norm-ultraweak hybrid
topology on matrix spaces, automatic strictness of certain
completely bounded maps,
 connections with the enveloping von Neumann algebra and ultraweak extension, and
 compatibility between extensions of maps continuous with respect to
these topologies.
 Here also basic properties of multiplier
$C^*$-bialgebras are described and a universal enveloping
construction is given.

The second part of the paper (Sections~\ref{Section:
QS},~\ref{Section: QSDE} and~\ref{Section: Cocycles}), is devoted to
developing the necessary quantum stochastic theory.
 Apart from its
immediate applications in this work, the contents of these sections
may be of wider interest to quantum stochastic analysts.
Section~\ref{Section: QS} contains a brief summary of the relevant
`standard' theory; weak and strong coalgebraic quantum stochastic
differential equations are treated in Section~\ref{Section: QSDE},
where an automatic strictness result is used to establish uniqueness
of weak solutions. Quantum stochastic convolution cocycles are
analysed in Section~\ref{Section: Cocycles}, where Markov-regular
completely positive contraction cocycles are shown to satisfy
quantum stochastic differential equations, and the form of the
stochastic generator is given, for these and for *-homomorphic
cocycles.

In the last part (Section~\ref{Section: Levy}), quantum L\'evy
processes are defined in our setting, and are shown to be realisable
as Fock space convolution cocycles when they are Markov-regular.
This leads to the characterisation of the generating functionals of
norm-continuous convolution semigroups of states on a locally
compact quantum semigroup $\Bil$ as functionals which have the form
\[
\gamma(b) = \langle \eta, \pi(b) \eta - \Cou(b) \eta \rangle, \;\; b \in \Bil,
\]
where $\pi$ is a nondegenerate representation of $\Bil$ on a Hilbert
space $\hil$ and $\eta$ is a vector in $\hil$.
 Examples of Poisson-type (i.e.\ norm-continuous) convolution
 semigroups of states on a locally compact quantum semigroup are
 thereby readily produced, and their associated quantum L\'evy
 processes constructed.

Our most satisfactory results are obtained in the case of
\emph{Markov-regular} quantum L\'evy processes.
A general theory of \emph{weakly continuous} convolution semigroups
of functionals on multiplier $C^*$-bialgebras is initiated
in~\cite{discrete}. In that paper every such semigroup of
\emph{states} on a multiplier $C^*$-bialgebra of discrete type is
shown to be automatically norm-continuous so that all the results of
this paper apply directly in that case. In~\cite{discrete} we have
used Theorem~\ref{significant} of this paper to derive a classical
result on conditionally positive-definite functions on a compact
group.

\emph{General notations.} In this paper the multiplier algebra and
universal enveloping von Neumann algebra of a $C^*$-algebra $\Al$
are denoted by $\Altilde$ and $\Alol$ respectively.
 The notations are intended to emphasise the view of these objects
 as two forms of completion of $\Al$, enjoying the natural
 inclusions $\Al \subset \Altilde \subset \Alol$.
 The symbols
$\otul$, $\ot$ and $\otol$ are used respectively for
linear/algebraic, spatial/minimal and ultraweak tensor products, of
spaces and respecively, linear, completely bounded and ultraweakly
continuous completely bounded maps. For any Hilbert space $\hil$, we
have the ampliation and Hilbert space, given respectively by
\begin{equation}
\label{iotahatDelta}
 \iota_\hil: B(\Hil; \Kil)\to B(\Hil\ot\hil;
\Kil\ot\hil), \ T\mapsto T\ot I_\hil, \ \text{ and } \hilhat:=
\Comp\oplus\hil.
\end{equation}
where context determines the Hilbert spaces $\Hil$ and $\Kil$.

\section{Strict extensions, tensor products and $\chi$-structure maps}
\label{Section: Extensions}
In this section we recall some definitions and relevant facts about
Hilbert $C^*$-modules (\cite{Lance}), strict topologies (\cite{Kus}),
tensor products and $\hil$-$\kil$-matrix spaces (\cite{LWexistence}).
We establish an automatic strictness result and show how strict
tensor product constructions compare with
$\hil$-matrix space constructions over a multiplier $\Cstar$-algebra.
The section ends by recalling a central concept for quantum L\'evy
processes, namely that of $\chi$-structure maps.

\subsection{Hilbert $\Cstar$-modules and multiplier
algebras}
For Hilbert $C^*$-modules $E$ and $F$ over a $C^*$-algebra $\Cil$,
$\Adjointable(E;F)$ denotes the space of adjointable operators $E\to F$.
Hilbert $\Cstar$-modules are endowed with a natural operator space structure
under which $M_n(\Adjointable(E;F))$ is identified with
$\Adjointable(E^n;F^n)$, where the \emph{column} direct sums $E^n$ and
$F^n$ are also Hilbert $\Cil$-modules, and
$\Adjointable(E;F) \subset CB_\Cil(E;F)$, the space of right $\Cil$-linear
completely bounded maps $E\to F$ (\cite{BlM}).
The \emph{strict topology} on $\Adjointable(E;F)$ is the locally convex
topology generated by the seminorms
$T\mapsto \|Te\| + \|T^*f\|$ ($e\in E, f\in F$);
it is Hausdorff and complete.
The closed subspace of
$\Adjointable(E;F)$ generated by the elementary maps
$| f\ra \la e|: x \mapsto f \la e, x\ra$
($e\in E$, $f\in F$) is denoted $\KAdjointable(E;F)$.
The unit ball of $\KAdjointable(E;F)$ is strictly dense in that of
$\Adjointable(E;F)$, $\KAdjointable(E)$ is a $\Cstar$-algebra and
$\Adjointable(E)$ is a model for its multiplier algebra. In
particular, viewing a $\Cstar$-algebra $\Al$ as a Hilbert
$\Cstar$-module over itself, $\KAdjointable(\Al) = \Al$ so that
$\Adjointable(\Al)$ is a model for the multipler algebra $\Altilde$.
A net of positive contractions $(e_i)_{i\in I}$ in $\Al$ is an
approximate identity for $\Al$ if and only if $e_i \to 1_{\Altilde}$
strictly. When $\Al$ is unital the strict topology coincides with
the norm topology. For Hilbert spaces $\hil$ and $\kil$, $|\hil\ra
:= B(\Comp;\hil)$ and $|\kil\ra$ are Hilbert $\Cstar$-modules over
$\Comp$, $\Adjointable(|\hil\ra ; |\kil\ra)$ and
$\KAdjointable(|\hil\ra ; |\kil\ra)$ are naturally identified with
$B(\hil;\kil)$ and $K(\hil;\kil)$ respectively, and the strict
topology on $\Adjointable(|\hil\ra ; |\kil\ra)$ corresponds to the
strong*-topology on $B(\hil;\kil)$. When a $\Cstar$-algebra $\Al$
acts nondegenerately on a Hilbert space $\hil$, the multiplier
algebra $\Altilde$ is realised as the double centraliser of $\Al$ in
$B(\hil)$: $\{x\in B(\hil): \forall_{a\in\Al} \ xa, ax \in \Al\}$,
the inclusion $\Altilde \subset \Al''$ holds and bounded strictly
convergent nets in $\Altilde$ converge (strong*- and thus)
$\sigma$-weakly.

Two elementary classes of strictly continuous maps that feature below
are component maps
$\varepsilon_{kl}:
\Adjointable (E_1\oplus E_2) \to \Adjointable(E_k;E_l)$,
$T = [T_{ij}] \mapsto T_{kl}$, where the column direct sum
$E_1\oplus E_2$ is a Hilbert $\Cil$-module, and multiplication operators
$\Adjointable(E;F) \to \Adjointable(E';F')$, $T\mapsto X^*TY$, where
$X\in \Adjointable(F';F)$ and $Y\in \Adjointable(E';E)$ for Hilbert
$\Cil$-modules $E'$ and $F'$.

\subsection{Strict maps and extensions}
There is a more prevalent notion in the theory than strict continuity. A
bounded operator $\varphi$ from
$\KAdjointable = \KAdjointable(E;F)$ to
$\Adjointable' = \Adjointable(E';F')$
where $E'$ and $F'$ are Hilbert $\Cstar$-modules over a $\Cstar$-algebra
$\Cil'$, is called \emph{strict} if it is strictly continuous on bounded
sets; the collection of such maps, denoted
$\Bstrict(\KAdjointable; \Adjointable')$, is a closed subspace of
$B(\KAdjointable; \Adjointable')$; we describe some of its contents below.
Here we are particularly interested in the classes $\Bstrict(\Al;
B)$ for $\Cstar$-algebras $\Al$ and spaces $B$ of the form $B(\hil;\kil)$
where $\hil$ and $\kil$ are Hilbert spaces.
An important general class of strict maps is the set of *-homomorphisms
$\varphi: \Al \to \Adjointable(E)$, for a $\Cstar$-algebra $\Al$ and
Hilbert $\Cstar$-module $E$,
which are \emph{nondegenerate} in the sense that
$\Linol\, \varphi(\Al)E = E$.

For a $C^*$-algebra $\Al$,
let $\Alol$ denote its universal enveloping von Neumann algebra,
let $\rho$ be the embedding $\Al \to \Alol$ and
let $\iota$ be the inclusion/natural map
$\Alol_* \to \Alol^* = (\Alol_*)^{**}$.
The map
$\iota^*\circ \rho^{**}: \Al^{**} \to (\Alol_*)^* = \Alol$
is a *-isomorphism for the common Arens product on $\Al^{**}$
and a weak*-$\sigma$-weak homeomorphism.
Since $\Al$ acts nondegenerately in the universal representation,
$\Altilde$ may be viewed
as a subalgebra of $\Alol$. All of this is well-known. For ease of
reference we collect together some extension properties which will play an
important role here.
The notation $\Buw$ stands for bounded ultraweakly continuous.

\begin{thm}
\label{Theorem 1.1}
Let $\Al$ be a $C^*$-algebra with multiplier algebra $\Altilde$
and universal enveloping von Neumann algebra $\Alol$.
\begin{alist}
\item
Let $\varphi \in \Bstrict(\Al;\Adjointable)$ where
$\Adjointable = \Adjointable(E;F)$ for $C^*$-modules $E$ and $F$ over a
$\Cstar$-algebra $\Cil$.
Then $\varphi$ has a unique strict
extension $\varphitilde: \Altilde \to \Adjointable$, moreover
$\varphitilde$ is bounded and $\|\varphitilde\| = \|\varphi\|$.
 \item
Let $\psi\in B(\Al;B)$ where $B = B(\hil;\kil)$ for Hilbert spaces $\hil$
and $\kil$. Then $\psi$ has a unique normal extension $\psiol \in
\Bsigma(\Alol; B)$, moreover $\|\psiol\|= \|\psi\|$.
\item
Let $\phi \in \Bstrict(\Al;B)$ where
$B=B(\hil;\kil)=\Adjointable(|\hil\ra;|\kil\ra)$ for
Hilbert spaces $\hil$ and $\kil$.
Then $\phitilde = \phiol|_{\Altilde}$.
\end{alist}
In \tu{(}a\tu{)},
$\widetilde{\varphi^\dagger} = \widetilde{\varphi}^\dagger$
where
$\varphi^\dagger: \Al \to \Adjointable(F;E)$ is defined by
$\varphi^\dagger(a^*) = \varphi(a)^*$; similarly, in \tu{(}b\tu{)}
$\ol{\psi^\dagger} = \psiol^\dagger$.
When $F=E$,
$\varphitilde$ is positive/completely positive/multiplicative if $\varphi$
is, and likewise for $\psi$ and $\psiol$ when $\kil = \hil$.
\end{thm}
\begin{proof}
(a) is proved in~\cite{Kus} in the case $E=F=\Cil$. The general case is
obtained by applying this case with $\Cil = \KAdjointable (E\oplus F)$ and
composing with the strict map $\Adjointable(E\oplus F) \to
\Adjointable(E;F)$, $T=[T_{i,j}]\mapsto T_{21}$.
(b) is well-known: set
$\psiol := \iota^* \circ \psi^{**}\circ j$ where $\iota$ is the natural
map/inclusion
$B_* \to (B_*)^{**} = B^*$ and $j$ is the natural isometric isomorphism
$\Alol \to \Al^{**}$.
Since the unit ball of $\Al$ is strictly dense in that of $\Altilde$
(\cite{Lance}, Proposition 1.4), (c)
follows from the fact that strictly convergent bounded nets converge
$\sigma$-weakly.
The last part follows from Kaplansky's density theorem
and its Hilbert $C^*$-module counterpart (just used),
and the separate strict
(respectively $\sigma$-weak) continuity of multiplication
and corresponding continuity of the adjoint operation, in a multiplier
algebra (respectively von Neumann algebra).
\end{proof}

\begin{rems}
(i)
The extensions commute with matrix liftings:
\begin{align*}
& \varphi^{(n)}\ \wt{ } =
 \varphitilde^{(n)}: M_n(\Al)\,\wt{ } =
M_n(\Altilde) \to M_n(\Adjointable)= \Adjointable(E^n;F^n)
\\ &
\ol{\psi^{(n)}} = \psiol^{(n)}: \ol{M_n(\Al)} = M_n(\Alol) \to
M_n(B)= B(\hil^n;\kil^n),
\end{align*}
so $\|\varphitilde\|_\cb = \|\varphi\|_\cb$ when
$\varphi\in \CBstrict(\Al;\Adjointable)$
and $\|\phiol\|_\cb = \|\phi\|_\cb$ when
$\psi\in CB(\Al;B)$.

(ii)
Clearly
the range of $\varphitilde$ is contained in the strict closure of
the range of $\varphi$, and
the range of $\psiol$ is contained in the $\sigma$-weak closure of
the range of $\psi$.

(iii)
As a consequence of (a), strict maps may be composed in the
following sense: if
$\varphi_1\in \Bstrict(\Al_1;\Altilde_2)$ and
$\varphi_2\in \Bstrict(\Al_2;\Altilde_3)$,
for $C^*$-algebras $\Al_1, \Al_2$ and $\Al_3$,
then
$a\mapsto \varphitilde_2 ( \varphi_1 (a))$
is strict
with unique strict extension $\varphitilde_2 \circ \varphitilde_1$;
following the widely adopted convention (e.g.\ \cite{Lance}), it is
simply denoted  $\varphi_2 \circ \varphi_1$, thus $\varphi_2 \circ
\varphi_1 \in \Bstrict(\Al_1;\Altilde_3)$.
 For a nondegenerate
*-homomorphism $\varphi: \Al \to \Ciltilde$, $\varphitilde$ is a
unital *-homomorphism, and conversely every nondegenerate
*-homomorphism $\varphi: \Al \to \Ciltilde$ is the restriction of a
strict unital *-homomorphism $\Altilde \to \Ciltilde$.
\end{rems}

\noindent \emph{Warning.} We now write
$\Bstrict(\Altilde;\Adjointable)$ for the class of strict maps
$\Altilde \to \Adjointable$ where, for us, $\Adjointable$ will
always be \emph{either} of the form $B(\hil;\kil) =
\Adjointable(|\hil\ra; |\kil\ra)$ for Hilbert spaces $\hil$ and
$\kil$ \emph{or} of the form $\Ciltilde$ for a $\Cstar$-algebra
$\Cil$ (or \emph{both}: $B(\hil) = \Adjointable(|\hil\ra) = K(\hil)\
\wt{ }$\ ). Use of this notation therefore always needs to reflect
the algebras of which the source (and target) multiplier algebras
are.

We note that the theorem delivers a commutative diagram of isometric
isomorphisms:
\begin{equation} \label{diag B}
\xymatrix{
\Bstrict(\Al;B) \ar [dr]\ar@{->}[rr] && \Buw(\Alol;B) \\
& \Bstrict(\Altilde;B) \ar[ur]&
}
\end{equation}
for any $\Cstar$-algebra $\Al$ and space $B$ of the form $B(\hil;\kil)$
for Hilbert spaces $\hil$ and $\kil$.

\begin{defn}
Let $\varphi\in\Bstrict(\Al; \Adjointable(E))$ for a $\Cstar$-algebra
$\Al$ and Hilbert $\Cstar$-module $E$. We call $\varphi$
\emph{preunital} if its strict extension is unital:
$
\varphitilde(1) = I.
$
\end{defn}
\begin{rem}
For *-homomorphisms this is equivalent to nondegeneracy; in general
it is equivalent to $\varphi(e_\lambda)\to 1$ for some/every
$\Cstar$-approximate identity $(e_i)_{i\in I}$ for $\Al$, but is
stronger than the condition $\Linol \varphi(\Al)E = E$
(\cite{Lance}, Proposition 2.5; Corollary 5.7).
\end{rem}

\subsection{Automatic strictness and strict tensor products}
In the next theorem we establish an automatic strictness property and
identify a natural class of maps for strict tensoring.
First some notation. When  $\Al$ is the spatial tensor product
$\Al_1 \ot \Al_2$, for $C^*$-algebras $\Al_1$ and $\Al_2$,
$\Altilde$ is denoted $\Al_1 \ottilde \Al_2$. Note the relation
\begin{equation} \label{ottilde}
\Altilde_1 \ot \Altilde_2 \subset \Al_1 \ottilde \Al_2
\end{equation}

\begin{thm} \label{P.B}
Let $\Al$, $\Al_1$ and $\Al_2$ be $C^*$-algebras.
\begin{alist}
\item
Let $\varphi \in CB(\Al;B)$ where $B=B(\hil;\kil)$ for a
Hilbert spaces $\hil$ and $\kil$. Then $\varphi$ is strict and
$\varphitilde = \varphiol|_{\Altilde}$.
In particular, all *-homomorphisms $\Al \to B(\hil)$ are strict.
\item
Let $\varphi_i\in \Lin \CPstrict\big(\Al_i; \Ciltilde_i\big)$ for
$C^*$-algebras $\Cil_1$ and $\Cil_2$.
 Then, there is a unique map
$\varphi_1\ot \varphi_2 \in \Lin \CPstrict(\Al_1\ot\Al_2;
\Cil_1\ottilde\Cil_2)$ extending the algebraic tensor product map
$\varphi_1 \otul \varphi_2$.
\end{alist}
\end{thm}

\begin{proof}
(a)
In view of Theorem~\ref{Theorem 1.1} it suffices to prove that $\varphi$
is
strict. It follows from the Wittstock-Paulsen-Haagerup Decomposition
Theorem (\cite{EfR}, Theorem 5.3.3)
that $\varphi = \psi \circ \pi$ where $\pi$ is a *-homomorphism
$\Al \to B(\Hil)$ for some Hilbert space $\Hil$ and
$\psi: B(\Hil) \to B$ is of the form $X\mapsto R^*XS$, for some
operators $R\in B(\kil;\Hil)$ and $S\in B(\hil;\Hil)$.
Moreover, replacing $\Hil$ by $\Hil':= \Linol \pi(\Al)\Hil$,
$\pi$ by its compression to $\Hil'$
$R$ and $S$ by $RP$ and $SP$
where $P$ is the orthogonal projection $\Hil \to \Hil'$, if necessary,
we may suppose that $\pi$ is
nondegenerate and therefore strict, when viewed as a map into
$\Adjointable(|\Hil\ra)$. Since $\psi$ is strict
$\Adjointable(|\Hil\ra) \to \Adjointable(|\hil\ra;|\kil\ra)$,
$\varphi$ is too.
(b)
By linearity we may suppose that $\varphi_1$ and $\varphi_2$ are
completely positive. The result then follows easily from Kasparov's
extension of the Stinespring Decomposition Theorem (\cite{Lance}, Theorem
5.6).
\end{proof}
\begin{rems}
It follows from Part (a) and the remarks after Theorem~\ref{Theorem 1.1}
that~\eqref{diag B} restricts to a commutative diagram of completely
isometric isomorphisms
\begin{equation} \label{diag CB}
\xymatrix{
CB(\Al;B) \ar [dr]\ar@{->}[rr] && \CBuw(\Alol;B) \\
& \CBstrict(\Altilde;B) \ar[ur]&
}
\end{equation}
where $\Al$ and $B$ are as in~\eqref{diag B}. In particular,
we have complete isometries
\begin{equation}
\label{Equation: duals} \Al^* \cong \Altilde_\beta^* \cong \Alol_*. 
\end{equation}

By operator space considerations
$\Ran(\varphi_1\ot\varphi_2) \subset \Ciltilde_1\ot\Ciltilde_2$.
Part (b) leads to the following useful notation. For
$\varphi_i\in \Lin \CPstrict\big(\Al_i; \Ciltilde_i\big)$ ($i=1,2$) we
denote the unique strict extension of
$\varphi_1\ot \varphi_2$ by
$\varphi_1\ottilde \varphi_2$. Thus
\begin{equation} \label{tensor tilde}
\varphi_1\ottilde \varphi_2 \in \Lin \CPstrict
\big(\Al_1\ottilde\Al_2;\Cil_1\ottilde\Cil_2\big).
\end{equation}
\end{rems}
Note the following consequence of
Part (a) and its proof,
which provides a source for tensoring as in Part (b).

\begin{cor}
\label{cor: CB=LinCPbeta}
For a $C^*$-algebra $\Al$ and Hilbert space $\hil$,
\begin{equation}
\label{CB=LinCPbeta}
CB\big(\Al; B(\hil)\big) = \Lin \, \CPstrict\big(\Al; B(\hil)\big).
\end{equation}
\end{cor}

\begin{rem}
This `strict decomposability' property is very useful. For  a
general multiplier algebra target space, completely bounded maps
need neither be strict nor be linear combinations of completely
positive maps.
\end{rem}

\subsection{$\hil$-$\kil$-Matrix spaces}
Let $\Vil$ be an operator space in $B(\Hil;\Kil)$, set $B =
B(\hil;\kil)$ for two further Hilbert spaces $\hil$ and $\kil$ with
total subsets $S$ and $T$, and let $Z\in B(\Hil\ot\hil;\Kil\ot\kil)
= B(\Hil;\Kil)\otol B$. Then the following are equivalent:
\begin{align*}
&
E^\xi Z E_\eta \in \Vil \text{ for all } \xi\in T, \eta \in S;
\\ &
(\id_{B(\Hil;\Kil)} \otol \omega)(Z) \in\Vil \text{ for all } \omega \in B_*;
\end{align*}
where, for a Hilbert space vector $\xi$,
\begin{equation}
\label{E notation}
E_{\xi} := I \ot |\xi\ra : u \mapsto u \ot \xi \text{ and }
E^{\xi} := (E_\xi)^* = I \ot \la\xi |,
\end{equation}
and $I$ is the identity operator on the appropriate Hilbert space. The collection of operators $Z$
enjoying this property is an operator space which is denoted $\Vil\otM B$ and called the (right)
$\hil$-$\kil$-\emph{matrix space over} $\Vil$. It is situated between the norm-spatial and
ultraweak-spatial tensor products:
\[
\Vil\ot B \subset
\Vil\otM B \subset
\ol{\Vil}\otol B
\]
and the latter inclusion is an equality if and only if $\Vil$ is
$\sigma$-weakly closed.
If $\varphi \in CB(\Vil;\Vil')$ for another concrete operator space
$\Vil'$
then there is a unique map, its $\hil$-$\kil$-\emph{lifting},
denoted $\varphi \otM \id_{B}$, from
$\Vil\otM B$ to $\Vil'\otM B$ satisfying
$E^\xi(\varphi\otM\id_B)(Z)E_\eta = \varphi(E^\xi ZE_\eta)$
for all $\xi\in\kil$, $\eta \in \hil$ and $Z \in \Vil\otM B$; it
is completely bounded, with $\|\varphi\otM\id_{B}\|_\cb =
\|\varphi\|_\cb$ and completely isometric if $\varphi$ is (unless
$B = \{0\}$), moreover it extends
$\varphi \ot \id_{B}$, and it coincides with
$\varphi \otol \id_{B}$ when $\Vil$ is $\sigma$-weakly
closed and $\varphi$ is $\sigma$-weakly continuous.
The next proposition confirms the compatibility of $\hil$-$\kil$-matrix spaces and
$\hil$-$\kil$-liftings on the one hand, and strict tensor products of
algebras and strict maps on the other.
First note the identity
\begin{equation}
\label{matrix ket id}
\|\eta\|^2 R E_\eta X =
R\big( X \ot |\eta \ra\la\eta|\big) E_\eta,
\end{equation}
for Hilbert space operators
$R\in B(\Hil\ot\hil;\Hil)$ and $X\in B(\Hil)$
and vectors $\eta\in\hil$.

\begin{propn}
Let $\Al$ be a $\Cstar$-algebra and let $\hil$ be a Hilbert space.
Then, in any faithful nondegenerate representation of $\Al$
\tu{(}such as its universal representation\tu{)},
\[
\Al\ottilde K(\hil) \subset \Altilde\otM B(\hil),
\]
for the induced concrete realisations of $\Altilde$ and $\Al\ottilde
K(\hil)$. Moreover, if $\psi\in\Lin\CPstrict(\Al;\Ciltilde)$ for
another $\Cstar$-algebra $\Cil$, also faithfully and nondegenerately
represented, then
\[
\psi \ottilde \id_{K(\hil)} \subset \psitilde \otM \id_{B(\hil)}.
\]
\end{propn}
\begin{proof}
Set $B=B(\hil)$ and $K=K(\hil)$. Let $T\in \Al\ottilde K$ and
$\zeta, \eta \in \hil$. First note that~\eqref{matrix ket id}
implies that
\begin{equation}
\label{to get strict}
\|\eta\|^2 E^\zeta T E_\eta a =
E^\zeta T \big( a \ot |\eta \ra\la\eta| \big) E_\eta,
\quad a \in \Al;
\end{equation}
similarly,
\begin{equation}
\label{2 to get strict}
\|\zeta\|^2 a E^\zeta T E_\eta  =
E^\zeta \big( a \ot |\zeta \ra\la\zeta| \big) T E_\eta,
\quad a \in \Al;
\end{equation}
and so $E^\zeta T E_\eta \in \Altilde$.
Thus $T\in\Altilde\otM B$.
This proves that $\Al \ottilde K \subset \Altilde\otM B$,
and~\eqref{to get strict} and~\eqref{to get strict}
now imply that the map
\[
T \in \Al \ottilde K \mapsto E^\zeta T E_\eta \in \Altilde
\]
is strictly continuous.
Therefore the maps
\[
E^\zeta \big( \psitilde \otM \id_B \big) (\cdot ) E_\eta =
\psitilde (E^\zeta \cdot E_\eta) \text{ and }
E^\zeta \big( \psi \ottilde \id_K \big) (\cdot ) E_\eta
\]
are strictly continuous $\Al \ottilde K \to \Ciltilde$
and agree on the strictly dense subspace $\Al \ot K$.
They therefore agree on $\Al \ottilde K$ and the result follows.
\end{proof}

%

\begin{rem}
In the universal representation of $\Al$ we have
the further compatibility relations,
\[
\Altilde\otM B(\hil) \subset \Alol \otol B(\hil)
\text{ and }
\psitilde \otM \id_{B(\hil)} \subset \psiol \otol \id_{B(\hil)}.
\]
\end{rem}

\subsection{$C^*$-algebras with character}
The following notion plays an important role in the theory.
Recall the notation~\eqref{iotahatDelta}.

\begin{defn}
A $\chi$-\emph{structure map on} a $\Cstar$-algebra with
character
$(\Al,\chi)$
is a linear map $\varphi : \Al \to B(\hilhat)$, for some Hilbert
space $\hil$, satisfying
\begin{equation}
\label{Equation: chi structure relation}
\varphi (a^*b) =
\varphi(a)^* \chi (b) + \chi(a)^* \varphi(b) +
\varphi (a)^* \Delta \varphi (b),
\end{equation}
where $\Delta:=
\left[\begin{smallmatrix} 0 & \\ & I_{\hil} \end{smallmatrix}\right]
 \in B(\hilhat)$ (no relation to coproducts).
\end{defn}

In terms of its block matrix form
 $\varphi = \left[\begin{smallmatrix} \gamma
& \mu \\ \lambda & \nu
\end{smallmatrix}\right]$, $\lambda$ is a kind of derivation (see
below) and $\mu = \lambda^\dagger: a \mapsto \lambda(a^*)^*$. More
specifically, the following result, established in~\cite{QSCC2},
gives the general form of $\chi$-structure maps.

\begin{thm}
 \label{established}
Let $(\Al,\chi)$ be a $\Cstar$-algebra with character and let
$\varphi$ be a linear map $\Al \to B(\hilhat)$, for some Hilbert
space $\hil$. Then the following are equivalent.
\begin{rlist}
\item
$\varphi$ is a $\chi$-structure map.
\item
$\varphi$  has block matrix form
\begin{equation}
\label{Equation: block} a \mapsto
\begin{bmatrix} \gamma(a) & \la \xi | \nu(a) \\
\nu(a) | \xi \ra & \nu(a) \end{bmatrix}
\text{ where } \gamma:= \omega_\xi \circ \nu \text{ for } \nu := \pi
- \iota_\hil \circ \chi,
\end{equation}
in which $(\pi, \hil)$ is a representation of $\Al$ \tu{(}as a
$\Cstar$-algebra\tu{)} and $\xi$ is a vector in $\hil$.
\end{rlist}
If $\varphi$ is a $\chi$-structure map with such a block matrix form
then it is necessarily strict, moreover $\pi$ is nondegenerate if
and only if $\varphitilde (1) = 0$.
\end{thm}

\begin{proof}
The first part is Theorem A6 of~\cite{QSCC2}.
 This implies that
$\varphi$ is completely bounded and so, by Theorem~\ref{P.B},
$\varphi$ is strict. After strict extension, the last part now
follows by inspection.
\end{proof}

We say that the $\chi$-structure map $\varphi$ is \emph{implemented}
by the pair $(\pi, \xi)$. Note the following alternative expression:
\[
 \varphi(a) =
\begin{bmatrix} \la \xi | \\ I_{\hil} \end{bmatrix}
\big( \pi(a) - \chi(a) I_\hil \big)
\begin{bmatrix} | \xi \ra & I_{\hil} \end{bmatrix}
 \quad (a \in \Al)
\]

\begin{rems}
Thus $\lambda$ is a $\pi$-$\chi$-derivation, in other words
$\lambda(ab) = \lambda(a)\chi(b) + \pi(a)\lambda(b)$ ($a,b\in\Al$),
which is implemented.

By the separate strict/$\sigma$-weak continuity of multiplication,
it follows that if $\varphi$ is a $\chi$-structure map then
$\varphitilde$ is a $\chitilde$-structure map and $\varphiol$ is a
$\chiol$-structure map.
\end{rems}

We shall need the following result in Section~\ref{Section: Levy}.
\begin{lemma} \label{condpos}
Let $(\Al,\chi)$ be a $C^*$-algebra with character.
Then, for any functional $\gamma \in \Al^*$,
if $\gamma$ is positive on $\Ker \chi$ then
$\wt{\gamma}$ is positive on $\Ker \chitilde$
and $\ol{\gamma}$ is positive on $\Ker \chiol$.
\end{lemma}

\begin{proof}
It suffices to prove that
$\Al_+ \cap \Ker \chi$ is strictly dense in $\Altilde_+ \cap \Ker \chitilde$
and
$\sigma$-weakly dense in $\Alol_+ \cap \Ker \chiol$.
Let $a \in \Altilde_+ \cap \Ker \chitilde$.
The Kaplansky Density Theorem for multiplier algebras (\cite{Lance}, Proposition 1.4)
implies that there is a bounded net $(c_i)_{i \in I}$
of selfadjoint elements in $\Al$ converging strictly to $a^{1/4}$.
Set $a_i = b_i^*b_i$ where
\[
b_i :=
c_i \big(c_i - \chi (c_i)\big) \in \Ker \chi.
\]
Then $a_i \in \Al_+ \cap \Ker \chi$ and separate strict continuity
of multiplication on bounded subsets of $\Al$, and strictness of
$\chi$, imply that $(a_i)_{i \in I}$ converges strictly to $a$. The
ultraweak density of $\Al_+ \cap \Ker \chi$ in $\Alol_+ \cap \Ker
\chiol$ is proved similarly, by appealing to the standard Kaplansky
Density Theorem (for von Neumann algebras).
\end{proof}

\section{Multiplier $C^*$-bialgebras}
\label{Section: Bialgebras}
It is convenient to consider bialgebras in both the $\Cstar$- and
$W^*$- categories and a universal enveloping operation linking the two.

\begin{defn}
A (\emph{multiplier}) \emph{$C^*$-bialgebra} is a
$C^*$-algebra $\Bil$ with \emph{coproduct},
that is a nondegenerate *-homomorphism
$\coproduct: \Bil \to \Bil\ottilde\Bil$ satisfying the coassociativity
conditions
\[
(\id_\Bil \ot \coproduct)\circ \coproduct =
(\coproduct \ot \id_\Bil)\circ \coproduct.
\]
A \emph{counit} for $(\Bil,\coproduct)$ is a character $\counit$ on $\Bil$
satisfying the counital property:
\[
(\id_\Bil \ot \counit)\circ \coproduct =
(\counit \ot \id_\Bil)\circ \coproduct = \id_\Bil.
\]
\end{defn}
\begin{rems}
The above definitions extend those for unital $C^*$-bialgebras, for which
$\Biltilde = \Bil$ and $\Bil\ottilde\Bil = \Bil\ot\Bil$.
The strict extension of a coproduct is a unital *-homomorphism and the
strict extension of a counit is a character on $\Biltilde$.
Note however that, in general, $(\Biltilde, \coproducttilde)$ is
\emph{not} itself a $\Cstar$-bialgebra as the inclusion
$\Biltilde \ot \Biltilde \subset \Bil \ottilde \Bil$
is usually proper.
\end{rems}

Examples of counital $\Cstar$-bialgebras include locally compact quantum
groups in the universal setting (\cite{Kus2}), in particular all
coamenable locally compact quantum groups are included.
If the assumptions on the coproduct $\coproduct$ are weakened to it being completely
positive, strict and preunital
then the resulting structure is called a (\emph{multiplier})
$\Cstar$-\emph{hyperbialgebra}
(cf.\ \cite{ChV}).

Let $\Bil$ be a $\Cstar$-bialgebra.
The \emph{convolute} of
$\phi_1\in\Lin \CPstrict(\Bil;\Altilde_1)$ and
$\phi_2\in\Lin \CPstrict(\Bil;\Altilde_2)$
for $\Cstar$-algebras $\Al_1$ and $\Al_2$
is defined by
\[
\phi_1\conv\phi_2 := (\phi_1\ot\phi_2)\circ \coproduct \in \Lin
\CPstrict(\Bil;\Al_1\ottilde\Al_2).
\]
We denote its strict extension by $\phi_1\convtilde\phi_2$.
Associativity of both of these convolutions follows from the
associativity of $\ottilde$ and coassociativity of
$\coproducttilde$. For each $C^*$-algebra $\Al$ define a map
\[
R_\Al: \Lin \CPstrict(\Bil;\Altilde) \to \CBstrict(\Bil;\Bil\ottilde\Al),
\quad
\phi \to
\id_\Bil \conv\phi = ( \id_\Bil \ot \phi ) \circ \coproduct.
\]
In case $\Al=\Comp$, $\Lin \CPstrict (\Bil;\Altilde)$ is simply $\Bil^*$
and we have
\[
R_{\Comp}(\varphi_1\conv\varphi_2) =
R_{\Comp}\varphi_1
\circ R_{\Comp}\varphi_2,
\quad \varphi_1, \varphi_2 \in \Bil^*.
\]
When $\Bil$ is counital each $R_\Al$ has left-inverse
\begin{equation}
\label{Equation: counit circ Rmap}
E_\Al: \CBstrict(\Bil;\Bil\ottilde\Al) \to \CBstrict(\Bil;\Altilde), \quad
\psi \to (\counit \ot \id_\Bil)\circ \psi.
\end{equation}

\begin{rems}
By the complete positivity and strictness of the coproduct
\[
R_\Al \big( \CPstrict (\Bil; \Altilde) \big) \subset
\CPstrict \big( \Bil; \Bil\ottilde \Al \big)
\]
for any $\Cstar$-algebra $\Al$.
In particular, by~\eqref{CB=LinCPbeta},
\[
R_{K(\hil)} \big( CB ( \Bil; B(\hil) ) \big) \subset
\Lin \CPstrict \big( \Bil; \Bil\ottilde K(\hil) \big).
\]
Note also that, when $\varphi_1 \in CB(\Bil;B(\hil_1))$ and
$\varphi_2 \in CB(\Bil;B(\hil_2))$
for Hilbert spaces $\hil_1$ and $\hil_2$,
\[
\varphi_1\conv\varphi_2 \in CB(\Bil;B(\hil_1\ot\hil_2)).
\]
\end{rems}

For convenience we summarise useful properties of the
$\mathcal{R}$-maps next.

\begin{propn}
\label{Proposition: R-map}
Let $\Bil$ be a $\Cstar$-bialgebra and let $\Al$ be a $\Cstar$-algebra.
Then $R_\Al$ is a completely contractive map into
$\CBstrict(\Bil; \Bil\ottilde\Al)$ with image in the subspace
$\Lin CP_\beta(\Bil; \Bil \ottilde \Al)$ and,
after strict extension, $R_{\Comp}$ is furthermore a
homomorphism of Banach algebras:
$(\Bil^*,\conv) \cong
\big((\Biltilde)^*_\beta, \convtilde \big) \to \CBstrict(\Biltilde)$.
When $\Bil$ is counital, $R_\Al$ is completely isometric with completely
contractive left-inverse $E_\Al$ and $R_{\Comp}$ is furthermore a unital algebra morphism.
\end{propn}

We now turn briefly to the $W^*$-category.
\begin{defn}
A \emph{von Neumann bialgebra} is a von Neumann algebra $\Mil$ with
coproduct, that is
a normal unital *-homomorphism $\coproduct: \Mil \to \Mil\otol\Mil$ which
is coassociative:
\[
(\id_\Mil\otol\coproduct)\circ \coproduct =
(\coproduct\otol \id_\Mil)\circ \coproduct.
\]
A counit for $(\Mil,\coproduct)$ is a normal character $\counit$ on $\Mil$
satisfying
\[
(\id_\Mil\otol\counit)\circ \coproduct =
(\counit\otol \id_\Mil)\circ \coproduct = \id_\Mil.
\]
\end{defn}
Convolution in this category is straightforward. Let
$\phi_1\in\CBuw(\Mil; Z_1)$ and $\phi_2\in\CBuw(\Mil; Z_2)$
for $\sigma$-weakly closed concrete operator spaces
$Z$, $Z_1$ and $Z_2$, then
\[
\phi_1\conv\phi_2 := (\phi_1\otol\phi_2)\circ \coproduct \in
\CBuw(\Mil; Z_1 \otol Z_2),
\]
so that we may define a map
\[
R^\sigma_Z: \CBuw( \Mil; Z ) \to \CBuw( \Mil; \Mil \otol Z ), \quad
\phi \mapsto \id_\Mil\conv\phi =
\big( \id_\Mil \otol \phi \big) \circ \coproduct.
\]
In particular,
\[
R^\sigma_\Comp(\varphi_1\conv\varphi_2) =
R^\sigma_\Comp\varphi_1 \circ R^\sigma_\Comp\varphi_2 \in CB_\sigma(\Mil),
\text{ for } \varphi_1, \varphi_2\in\Mil_*.
\]
When $\Mil$ is counital
$R^\sigma_Z$ has left-inverse
\[
E^\sigma_Z: \CBuw( \Mil; \Mil \otol Z ) \to \CBuw( \Mil; Z ), \quad
\psi \mapsto (\counit \otol \id_Z )\circ \psi.
\]
\begin{propn}
Let $(\Bil, \coproduct)$ be a $\Cstar$-bialgebra.
Then $(\Bilol, \coproductol)$ is a von Neumann bialgebra.
Moreover, if $\counit$ is a counit for $\Bil$ then $\counitol$ is a counit
for $\Bilol$.
\end{propn}
\begin{proof}
The map $\coproductol$ is a normal, unital *-homomorphism and the normal
maps
\[
(\id_{\Bilol}\otol\coproductol)\circ \coproductol \text{ and }
(\coproductol\otol \id_{\Bilol})\circ \coproductol.
\]
agree on $\Bil$, which is $\sigma$-weakly dense in $\Bilol$, and so
coincide. In the counital case, $\counitol$ is a normal character on
$\Bilol$ and the normal maps
\[
\big( \id_{\Bilol} \otol \counitol \big) \circ \coproductol, \
\big( \counitol\otol \id_{\Bilol} \big) \circ \coproductol
\text{ and } \id_{\Bilol}
\]
agree on $\Bil$ and so coincide.
\end{proof}
Naturally, we refer to $(\Bilol, \coproductol)$, respectively
$(\Bilol, \coproductol, \counitol)$
as the \emph{universal enveloping von Neumann bialgebra}
(resp. \emph{counital von Neumann bialgebra}) \emph{of} $\Bil$.
\begin{rem}
The two forms of $\mathcal{R}$-map enjoy an easy compatibility:
if $\phi \in \Lin CP_\beta\big(\Bil;\Altilde\big)$ for a $\Cstar$-algebra $\Al$
then
\begin{equation}
\label{compatibility of R maps}
\ol{R_\Al \phi} = R^\sigma_{\Alol} \,\phiol,
\end{equation}
and similarly for the $E$ maps in the counital case.
\end{rem}

\emph{From now on} we shall denote all maps of the form
$R^\sigma_Z$, respectively $E^\sigma_Z$, by
$\Rbar$, respectively $\Ebar$,
and similarly abbreviate all maps of the form $R_\Al$ and $E_\Al$ to
$\Rmap$ and $\Emap$.

\section{Quantum stochastics}
\label{Section: QS}

\emph{Fix now, and for the rest of the paper},
a complex Hilbert space $\kil$ referred to as the
\emph{noise dimension space}.
For $c\in \kil$ define $\chat:=\binom{1}{c}\in \kilhat$;
and for any function $g$ with values in $\kil$
let $\ghat$ denote the corresponding
function with values in $\kilhat$, defined by $\ghat(s):= \widehat{g(s)}$.
Let $\FFock$ denote the symmetric Fock space over $L^2(\Rplus;\kil)$, let
$\Step$ denote the linear span of $\{d_{[0,t[}: d\in \kil, t\in\Rplus\}$
in $L^2(\Rplus;\kil)$ (for purposes of evaluating, we always take these
right-continuous versions) and
let $\Exps$ denote the linear span of $\{\ve(g): g\in\Step\}$ in $\FFock$,
where $\ve(g)$ denotes the exponential vector
$\big((n!)^{-\frac{1}{2}}g^{\ot n}\big)_{n\geq 0}$.
(There will be no danger of confusion with the inverse of an
$\mathcal{R}$-map!) Also define
\begin{equation}
\label{Delta QS}
e_0 :=
\binom{1}{0}\in\kilhat
\text{ and }
\dQS := P_{\{0\}\oplus\kil} =
\begin{bmatrix} 0 & \\ & I_{\kil}
\end{bmatrix} \in B(\kilhat).
\end{equation}

\subsection{Quantum stochastic processes, differential equations and cocycles}
A detailed summary of the relevant results from QS analysis
([$\text{LW}_{\!1\text{-}4}$],~\cite{LSqsde}) is given
in~\cite{QSCC2}. We shall therefore be brief here.

For operator spaces $\Vil$ and $\Wil$, with $\Wil$ concrete,
$\Proc(\Vil \to \Wil)$ denotes the space of adapted proceses $k =
(k_t)_{t\geq 0}$ thus, for $t\in\Rplus$,
\begin{equation}
 \label{kte}
 k_t \in L\big(\Exps; L(\Vil; \Wil\otM |\FFock\ra)\big),
 \text{ written }
  \varepsilon\mapsto k_{t,\varepsilon}.
\end{equation}
 As in~\cite{QSCC2}, we abbreviate to $\Procstar(\Vil)$ when
$\Wil=\Comp$. Its associated maps $\kappa^{f,g}_t: \Vil\to\Wil$
($f,g \in \Step, s,t \in \Rplus$) are defined by
\[
\kappa^{f',f}_t(x)  =
\big(\id_\Wil \otM \la\varepsilon'|\big) k_{t,\varepsilon}(x),
\quad x\in\Vil,
\]
where $\varepsilon = \varepsilon(f_{[0,t[})$ and
$\varepsilon' = \varepsilon(f'_{[0,t[})$.
For us here
the pair $(\Vil,\Wil)$ will be either $(\Mil,\Comp)$ or $(\Mil,\Mil)$
for a von Neumann algebra $\Mil$, or $(\Al,\Comp)$ or
$(\Altilde,\Comp)$ for a $C^*$-algebra $\Al$ with multiplier algebra
$\Altilde$. Thus $\Wil \otM |\Fock\ra$ is either $|\Fock\ra$ or
$\Mil \otol |\Fock\ra$.
The process $k$ is \emph{weakly initial space bounded} if each
$\kappa^{f,g}_t$ is bounded, and
\emph{weakly regular} if further
$\sup \big\{ \| \kappa^{f,g}_s  \| : s\in [0,T] \big\} < \infty$,
for all $T \geq 0$. Here `column-boundedness' usually obtains:
$k_{t,\varepsilon} \in B(\Vil; \Wil\otM |\Fock\ra)$ or
$CB(\Vil; \Wil\otM |\Fock\ra)$,
for each $t\in\Rplus$ and $\varepsilon\in\Exps$.
For von Neumann algebras $\Mil$ and $\Nil$,
a process $k \in \Proc(\Mil\to\Nil)$
is called \emph{normal} if each
$\kappa^{f,g}_t$ is normal.
It follows that if $k\in\Proc(\Mil\to\Nil)$ is bounded (meaning that each
$k_t$ is bounded) and normal then each map $k_t$ is normal
$\Mil \to \Nil\otol B(\Fock)$.
Note that, by Theorem~\ref{P.B},
any completely bounded process $l \in \Procstar(\Al)$ on a $C^*$-algebra $\Al$
is necessarily strict in the sense that each map $l_t: \Al\to B(\FFock)$ is strict.

For $\phi\in L( \kilhat; CB(\Vil; \Vil\otM |\kilhat\ra)$
and $\kappa \in CB(\Vil;\Wil)$, where both $\Vil$ and $\Wil$ are concrete
operator spaces, $k^{\phi,\kappa}$ denotes the unique weakly regular
process $k\in\Proc(\Vil\to\Wil)$ satisfying the quantum stochastic
differential equation
\begin{equation}
\label{QSDE}
dk_t = k_t \fullcomp d\Lambda_\phi (t), \quad
k_0 = \iota_{\FFock} \circ \kappa ,
\end{equation}
where $\iota_\FFock$ denotes the ampliation $\Wil\to\Wil\otM B(\FFock)$
$x\mapsto x\ot I_\FFock$.
The solution is given by
\[
k_{t,\ve} =
\sum_{n\geq 0} \Lambda^n_{t,\ve} \circ (\kappa \circ \phi_n)
\]
where $\phi_n$ is an $n$-fold composition of matrix liftings of
$\phi$,
 $\Lambda^n_{t,\ve}: \Wil \otM
|B(\khat^{\ot n})\ra \to \Wil \otM |\FFock\ra$ is the $\ve$-column
of the $n$-fold multiple QS integration map,
 and the sum is norm-convergent in $CB(\Vil; \Wil \otM
|\FFock\ra)$.

When $\Wil = \Vil$ and $\kappa = \id_\Vil$, $k$ is a weak quantum
stochastic cocycle on $\Vil$ (denoted $k^\phi$), that is it satisfies
$k_0 = \iota_\FFock$ and
for $s,t\in\Rplus$ and $f,g\in\Step$,
\begin{equation}
\label{Equation: weak standard cocycle}
\kappa^{f,g}_0 = \id_\Vil, \quad
\kappa^{f,g}_{s+t} =
\kappa^{f,g}_s \circ
\kappa^{S^*_sf,S^*_sg}_t \quad
(f,g \in \Step, s,t \in \Rplus)
\end{equation}
where $(S_t)_{t\geq 0}$
is the isometric semigroup of right-shifts on $L^2(\Rplus; \kil)$.
Let $(\sigma_t)_{t\geq 0}$ denote the induced endomorphism
semigroup on $B(\FFock)$, ampliated to $\Wil \otM B(\FFock)$.
Then, when $k$ is a completely bounded process, the cocycle relation
simpiflies to
\[
k_{s+t} = k_s \fullcomp \sigma_s \circ k_t,
\]
where the extended composition notation (which we do not need to go into
here) is explained in~\cite{QSCC2}.

\section{Coalgebraic quantum stochastic differential equations}
\label{Section: QSDE}

\emph{For this section we fix a $C^*$-bialgebra} $\Bil$,
which we do not assume to be counital, and
consider the coalgebraic quantum stochastic differential equation
\begin{equation}
\label{Equation: coalg QSDE}
dl_t = l_t \conv d\Lambda_\varphi (t), \quad
l_0 = \iota_{\FFock} \circ \eta ,
\end{equation}
where $\varphi\in SL(\khat, \khat; \Bil^*)$ and $\eta \in \Bil^*$.

\begin{defn}
By a \emph{form solution} of~\eqref{Equation: coalg QSDE}
is meant a family
$\big\{ \lambda^{f,g}_t \big| f, g \in \Step ,t \in \Rplus\big\}$
in $\Bil^*$ satisfying
\begin{rlist}
\item
the map
$s \mapsto \big( \lambda^{f,g}_s \conv \varphi_{\fhat(s),\ghat(s)} \big) (b)$
\ is locally integrable;
\label{i}
\item
$\lambda^{f,g}_t (b) - e^{\la f,g\ra} \eta (b) =
\int^t_0 \, ds \,
\big( \lambda^{f,g}_s \conv \varphi_{\fhat(s),\ghat(s)} \big) (b)$
\label{ii}
\end{rlist}
for all $f,g \in \Step$, $t \in \Rplus$ and $b \in \Bil$.
\end{defn}
\begin{rems}
Let $f,g \in \Step$ and $b \in \Bil$.
By automatic strictness of bounded linear functionals on $\Bil$,
(\ref{i}) makes sense.
By (\ref{ii})
it follows that $\lambda^{f,g}_t(b)$ is continuous in $t$,
and so is locally bounded. Therefore, by the Banach-Steinhaus Theorem,
$\lambda^{f,g}_t$ is locally bounded in $t$ and
(\ref{ii}) therefore implies that (\ref{i}) refines to
\begin{itemize}
\item[(i)$'$]
the map $s \mapsto \lambda^{f,g}_s$ is continuous,
\end{itemize}
which in turn implies that (\ref{ii}) refines to
\begin{itemize}
\item[(ii)$'$]
$\lambda^{f,g}_t  - e^{\la f,g\ra} \eta  =
\int^t_0 \, ds \,
\lambda^{f,g}_s \conv \varphi_{\fhat(s),\ghat(s)}$,
\end{itemize}
the integrand being piecewise norm-continuous $\Real_+ \to \Bil^*$.
\end{rems}

The following automatic strictness property is needed to establish
uniqueness for form solutions. Recall the strict extension notation
$\convtilde$.

\begin{lemma}
\label{lemma 5.1}
Let $\varphi\in SL(\khat, \khat; \Bil^*)$ and $\eta \in \Bil^*$.
Then every form solution
$\big\{ \lambda^{f,g}_t \big| f, g \in \Step ,t \in \Rplus\big\}$
of~\eqref{Equation: coalg QSDE}
is \emph{strict} in the sense that it satisfies
\begin{itemize}
\item[(ii)$\,\wt{}\ $]
$\ \widetilde{\lambda}^{f,g}_t  - e^{\la f,g\ra} \widetilde{\eta}  =
\int^t_0 \, ds \,
\lambda^{f,g}_s \convtilde \varphi_{\fhat(s),\ghat(s)}$
for all $f,g \in \Step$ and $t \in \Rplus$,
\end{itemize}
where
$\widetilde{\lambda}^{f,g}_t := \big( \lambda^{f,g}_t\big)\,\wt{}$.
\end{lemma}

\noindent
Note that the integrand in (ii)$\,\tilde{}$ is piecewise continuous
$\Rplus \to (\Biltilde)^*_\beta$.
\begin{proof}
Let
$\big\{ \lambda^{f,g}_t \big| f, g \in \Step ,t \in \Rplus\big\}$
be a form solution of~\eqref{Equation: coalg
QSDE}
and let $f,g \in \Step$ and $t \in \Rplus$.
Define bounded linear functionals
\[
\Phi :=
\int^t_0 \, ds \,
\lambda^{f,g}_s \conv \varphi_{\fhat(s),\ghat(s)} \text{ and }
\Psi :=
\int^t_0 \, ds \,
\lambda^{f,g}_s \convtilde \varphi_{\fhat(s),\ghat(s)}
\]
on $\Bil$ and $\Biltilde$ respectively. Note that each Riemann
approximant $\Psi_\mathcal{P}$ of $\Psi$ equals
$(\Phi_\mathcal{P})\,\wt{}$ where
$\Phi_\mathcal{P}$ is the corresponding Riemann approximant of $\Phi$.
The extension map $\Bil^* \to \Biltilde^*$ is (isometric and thus)
continuous therefore
\[
\Psi = \lim \Psi_\mathcal{P} =
(\lim \Psi_\mathcal{P})\,\widetilde{} = \Phi\,\widetilde{}.
\]
Since
$\Phi = \lambda^{f,g}_t - e^{\la f, g\ra} \eta$
it follows that
$\Psi = (\lambda^{f,g}_t)\,\widetilde{} - e^{\la f, g\ra}
\widetilde{\eta}$.
Thus
form solution is strict.
\end{proof}

With this we have uniqueness as well as existence for
form solutions.

\begin{thm}
\label{Theorem: SLQSDE}
Let $\varphi\in SL(\khat, \khat; \Bil^*)$ and $\eta \in \Bil^*$, for a
$C^*$-bialgebra $\Bil$. Then the quantum stochastic differential
equation~\eqref{Equation: coalg QSDE}
has a unique form solution.
\end{thm}
\begin{proof}
For each $c,d\in\kil$ let $(p^{c,d}_t)_{t\geq 0}$ denote the
norm-continuous one-parameter semigroup generated by
$\varphi_{\chat, \dhat} \in \Bil^*$, in the unitisation of the
Banach algebra $(\Bil^*, \conv)$. For $f,g\in\Step$ and
$t\in\Rplus$, set
\[
\lambda^{f,g}_t := \eta \conv p^{c_0,d_0}_{t_1-t_0} \conv  \cdots  \conv
p^{c_{n},d_n}_{t_{n+1}-t_n}
\]
where $t_0 = 0$, $t_{n+1} = t$, $\{t_1 < \cdots < t_n\}$ is the
(possibly empty) set of points in $]0,t[$ where $f$ or $g$ is
discontinuous and $(c_i,d_i) = (f(t_i),g(t_i))$ for $i=0, \cdots ,
n$. This defines an element $\lambda^{f,g}_t$ of $\Bil^*$. It is
easily verified that the resulting family $\big\{ \lambda^{f,g}_t
\big| f, g \in \Step ,t \in \Rplus\big\}$ is a form solution
of~\eqref{Equation: coalg QSDE}.

Suppose now that $\mu$ is the difference of two form
solutions, and let $f,g\in\Step$ and $t\in\Rplus$.
Then Lemma~\ref{lemma 5.1} yields
the identity
\[
\big(\mu^{f,g}_t\big)\,\wt{} =
\int^t_0 \, ds \,
\mu^{f,g}_s \convtilde \varphi_{\fhat(s),\ghat(s)},
\]
which may be iterated. Estimating after repeated iteration
(and using the isometry $\Biltilde^*_\beta \cong \Bil^*$) we have
\[
\|\mu^{f,g}_t\| \leq
\frac{t^n}{n!} \sup_{s\in[0,t]} \|\mu^{f,g}_s\|
\max \big\{ \| \varphi_{\chat,\dhat} \| :c \in \Ran f , d \in \Ran g
\big\}^n
\]
which tends to $0$ as $n \to \infty$.
Thus $\mu = 0$, proving uniqueness.
\end{proof}

We now show how stronger forms of solution are obtained when the
coefficient of the quantum stochastic differential equation is a
bounded mapping rather than just a form.
Below the following natural inclusions are invoked:
\begin{align*}
&B\big(\Bil; B(\hil; \hil')\big) \cong
B\big( \ol{\hil'}, \hil; \Bil^* \big) \subset
SL\big(\hil', \hil; \Bil^*\big),
\\ &
\varphi \mapsto
\big( (\zeta,\eta) \mapsto \varphi_{\zeta,\eta} :=
\la \zeta, \varphi(\cdot)\eta \ra \big),
\end{align*}
for Hilbert spaces $\hil$ and $\hil'$. Recall the notation for the
solution of a QS differential equation introduced above
equation~\eqref{QSDE}, and the column notation for
processes~\eqref{kte}.

\begin{thm}
\label{Theorem: lphieta}
Let $\varphi \in CB \big( \Bil ; B(\khat) \big)$
and $\eta \in \Bil^*$, for a $C^*$-bialgebra $\Bil$. Set
\[
\widetilde{l}^{\varphi, \eta} :=
k^{\widetilde{\phi}, \widetilde{\eta}} \text{ and }
\ol{l}^{\varphi, \eta} := k^{\ol{\phi}, \ol{\eta}}
\text{ where } \phi :=
\Rmap \varphi.
\]
Thus $\widetilde{\phi}
\in \CBstrict\big(\Biltilde; \Bil \ottilde K(\khat)\big)$ and
$\ol{\phi}
=
\Rbar
\ol{\varphi}
\in \CBuw\big(\Bilol; \Bilol \otol B(\khat)\big)$.
\begin{alist}
\item
Abbreviating
$\ol{l}^{\varphi, \eta}$ to $\ol{l}$ and
$\wt{l}^{\varphi, \eta}$ to $\wt{l}$
we have,
for all $\varepsilon \in \Exps$ and $t\in\Rplus$,
\begin{rlist}
\item
$\ol{l}_{t,\varepsilon} \in
\CBuw\big(\Bilol; |\FFock\ra\big)$\tu{;}
\item
$\wt{l}_{t,\varepsilon} =
\ol{l}_{t,\varepsilon}|_{\Biltilde}$\tu{;}
\item
$\wt{l}_{t,\varepsilon} \in
\CBstrict\big(\Biltilde; |\FFock\ra\big)$.
\end{rlist}
\item
For all $f,g\in\Step$ and $t\in\Rplus$,
setting
\begin{equation}
\label{Equation: lambda and kappa}
\ol{\lambda}^{f,g}_t :=
\omega_{\varepsilon(f_{[0,t[}), \varepsilon(g_{[0,t[})}
\circ \ol{l}^{\varphi, \eta}_t
\text{ and }
\kappa^{f,g}_t :=
E^{\varepsilon(f_{[0,t[})} k_t^{\ol{\phi}}( \cdot )
E_{\varepsilon(g_{[0,t[})},
\end{equation}
\begin{rlist}
\item
$\big\{ \ol{\lambda}^{f,g}_t|_{\Bil}\, \big|\ f, g \in \Step ,t \in \Rplus\big\}$
is the unique form solution of~\eqref{Equation: coalg QSDE};
\item
$\ol{\lambda}^{f,g}_t = \ol{\eta} \circ \kappa^{f,g}_t$ and
$
\Rbar
\,\ol{\lambda}^{f,g}_t =
(\Rbar
\,\ol{\eta}) \circ \kappa^{f,g}_t$.
\end{rlist}
\end{alist}
\end{thm}

\begin{proof}
Fix $\ve\in\Exps$ and $t\geq 0$. By linearity we may assume that
$\varepsilon = \varepsilon(g)$. Below we adopt the normal extension
notation $\alpha \convol \beta$ := $\ol{\alpha \conv \beta}$.

(a)
(i)
The operator
$\ol{l}_{t,\varepsilon}$ is a norm-convergent sum,
in $CB\big(\Bilol; |\FFock\ra\big)$,
of terms of the form
$\Lambda^n_{t,\ve} \circ \big(\eta\convol \varphi^{\convol n}\big)$
($n\in\ZZ_+$),
and each map $\eta\convol \varphi^{\convol n}$ is $\sigma$-weakly
continuous.
Since $\CBuw\big(\Bilol; |\FFock\ra\big)$ is a norm-closed subspace of
$CB\big(\Bilol; |\FFock\ra\big)$ it remains only to show that the bounded
operator
$\Lambda^n_{t,\ve}:
B\big(\khat^{\ot n}\big) \to |\FFock\ra$
is $\sigma$-weakly continuous. By the Krein-Smulian Theorem
it suffices to prove this on bounded sets.
This follows
from the following identity for multiple QS integrals:
\[
\big\la \varepsilon(f), \Lambda^n_t(A)\varepsilon(g) \big\ra =
\int_{\Delta^n_t} d\mathbf{s}
\, \big\la \pi_{\fhat}(\mathbf{s}), A \pi_{\ghat}(\mathbf{s}) \big\ra
e^{\la f, g \ra}, \quad A\in B(\khat^{\ot n}),
\]
since the integrand is a step function on
$\Delta^n_t :=
\{ s \in \mathbb{R}^n : 0 \leq s_1 \leq \cdots \leq s_n \leq t \}$.

(ii)
Since
$\wt{l}_{t,\varepsilon}$
is a norm-convergent sum, in $CB(\Biltilde; |\FFock\ra)$,
of terms of the form
$\Lambda^n_{t,\ve} \circ (\eta\convtilde \varphi^{\convtilde n})$,
this follows from (i) and the identity
$
\eta\convol \varphi^{\convol n}|_{\Biltilde} =
\eta\convtilde \varphi^{\convtilde n}
$
($n\in\ZZ_+$).

(iii)
This follows from (i) and (ii) since, for any map
$\alpha \in \CBuw( \Bilol; | \FFock \ra )$,
$(\alpha|_{\Bil})\ \wt{}
= (\alpha|_{\Biltilde})\ \wt{ }$\ (see~\eqref{diag CB}).

(b)
(i)
This follows from the identity
\[
\omega_{\varepsilon, \varepsilon'}\circ k^{\ol{\phi}, \ol{\eta}}_s
\circ \ol{\phi}_{\chat,\dhat}
=
\omega_{\varepsilon, \varepsilon'}\circ k^{\ol{\phi}, \ol{\eta}}_s
\, \convol \, \ol{\varphi}_{\chat,\dhat},
\quad
\varepsilon, \varepsilon'\in\Exps, c,d\in\kil, s\in\Rplus,
\]
where
$\ol{\phi}_{\chat,\dhat}:=
\big(\id_{\Bilol}\otol \omega_{\chat,\dhat}\big)\circ \ol{\phi}$.

(ii)
The first identity expresses the general relation between
$k^{\ol{\phi}, \ol{\eta}}$ and $k^{\ol{\phi}}$
(\cite{LSqsde}).
By (i), it follows from the proof of Theorem~\ref{Theorem: SLQSDE}
that
$\ol{\lambda}^{f,g}_t$ may be written in the form
\[
\ol{\eta} \, \convol \, \ol{p}^{c_0,d_0}_{t_1-t_0} \convol  \cdots
\convol \ol{p}^{c_{n},d_n}_{t_{n+1}-t_n},
\]
where $\ol{p}^{c,d}_t$ denotes the normal extension of $p^{c,d}_t$.
Thus
$
\Rbar
\ol{\lambda}^{f,g}_t$ equals
\[
\Rbar
\ol{\eta} \circ \ol{P}^{c_0,d_0}_{t_1-t_0} \circ  \cdots
\circ \ol{P}^{c_{n},d_n}_{t_{n+1}-t_n},
\]
where
$\ol{P}^{c,d}_t := \exp t\big(\phiol_{\chat,\dhat}\big) =
\Rbar
\ol{p}^{c,d}_t$.
(ii)
therefore now follows from the semigroup representation of the standard
QS cocycle $k^{\ol{\phi}}$ (\cite{LWjfa}).
\end{proof}

\begin{notn}
Setting $l^{\varphi ,\eta}_{t, \varepsilon} =
\ol{l}^{\varphi ,\eta}_{t, \varepsilon}|_{\Bil}$
($t\in\Rplus, \varepsilon\in\Exps$)
defines a process
$l^{\varphi ,\eta} \in \Procstar(\Bil)$,
which we denote by $l^{\varphi}$ when $\Bil$ is counital and
$\eta = \counit$.
This extends the notation introduced in~\cite{QSCC2} for the unital case.
\end{notn}

\begin{rem}
In view of the identity
\[
\big( k_s^{\phiol,\ol{\eta}}\circ \phi \big)_{\ve(f), \fhat(s)} =
\ol{l}_{s,\ve(f)}^{\varphi,\eta} \convol \varphiol_{\fhat(s)},
\]
$\ol{l}^{\varphi,\eta}$ satisfies
\[
\ol{l}_t = \iota_{\FFock} \circ \ol{\eta} +
\int_0^t \ol{l}_s \convol \varphiol \, d\Lambda(s), \quad t\in\Real_+.
\]
In this sense, $l^{\varphi,\eta}$ is a
strong solution of~\eqref{Equation: coalg QSDE}.


\end{rem}

Note that only the coalgebraic structure of $\Bil$ has been used so
far, not its algebraic structure.

We end this section by noting some correspondence between convolution
processes and associated standard processes.
Recall the notation for QSDE solutions
introduced above equation~\eqref{Equation: weak standard cocycle}.

\begin{propn}
\label{Proposition: XX}
Let $l = l^\varphi$ and $\ol{l} = \ol{l}^\varphi$,
where $\varphi\in CB(\Bil; B(\khat))$
for a counital $\Cstar$-bialgebra $\Bil$, and set
$k = k^{\ol{\phi}}$ where $\phi := \Rmap \varphi$.
Then
\begin{alist}
\item
$\ol{l}$ is unital if and only if $k$ is.
\item
$l$ is completely bounded \tu{(}respectively, completely positive
or *-homomorphic\tu{)} if and only if $k$ is, in which
case
\[
k_t =
\Rbar
\ol{l}_t, \quad
\ol{l}_t =
\Ebar
k_t \text{ and }
\|l_t\|_\cb = \|k_t\|_\cb, \quad t\in\Rplus.
\]
\end{alist}
\end{propn}

\begin{proof}
In the notations~\eqref{Equation: lambda and kappa},
Theorem~\ref{Theorem: lphieta}(b)(ii) implies that,
\[
\ol{\lambda}^{f,g}_t =
\Ebar
\kappa^{f,g}_t \text{ and }
\kappa^{f,g}_t =
\Rbar
\ol{\lambda}^{f,g}_t,
\quad f,g\in\Step, t\in\Rplus.
\]
Thus (a) follows from the unitality of the maps
$\counitol$ and $\coproductol$.
Moreover, if $k$ is completely bounded then, since
\[
\omega_{\ve,\ve'}
\circ \ol{l}^{\varphi, \eta}_t =
\ol{\lambda}^{f,g}_t =
\Ebar
\kappa^{f,g}_t =
\omega_{\ve,\ve'}
\circ
\Ebar
k_t,
\]
where $\ve = \varepsilon(f_{[0,t[}$ and $\ve' = \varepsilon(f'_{[0,t[}$,
for all $f,f'\in\Step$ and $t\in\Rplus$, it follows that
$l_t =
\Ebar
k_t$
($t\in\Rplus$), in particular $l$ is completely bounded.
Conversely, if $l$ is completely bounded then
$\ol{l_t} = \ol{l}_t$ ($t\in\Rplus$) and
\[
\big( \id_{\Bilol}\otol\omega_{\ve,\ve'} \big) \circ
\Rbar
\ol{l}_t
=
\Rbar
\ol{\lambda}^{f,f'}_t =
\kappa^{f,f'}_t =
\big( \id_{\Bilol}\otol\omega_{\ve,\ve'} \big) \circ
k_t
\]
for all $f,f'\in\Step$ and $t\in\Rplus$, so
$k_t =
\Rbar
\ol{l}_t$ ($t\in\Rplus$), therefore
$k$ is completely bounded.
The rest follows from the fact that $\coproductol$ and
$\counitol \otol \id_{B(\FFock)}$ are *-homomorphisms.
\end{proof}

\section{Quantum stochastic convolution cocycles}
\label{Section: Cocycles}

For this section \emph{we fix a counital $C^*$-bialgebra} $\Bil$.

\begin{defn}
A family $ \big\{ \lambda^{f,g}_t \big| f, g \in \Step ,t \in
\Rplus\big\}$ in $\Bil^*$ is a \emph{form quantum stochastic
convolution cocycle} on $\Bil$ if it satisfies
\[
\lambda^{f,g}_0 = \counit, \quad
\lambda^{f,g}_{s+t} =
\lambda^{f,g}_s \conv
\lambda^{S^*_sf,S^*_sg}_t, \quad
f,g \in \Step, s,t \in \Rplus,
\]
where
$(S_t)_{t\geq 0}$
is the isometric shift semigroup on $L^2(\Rplus; \kil)$.
\end{defn}
Note that for such a cocycle
\[
p_t^{c,d}:= \lambda_t^{c_{[0,t[}, d_{[0,t[}}
\]
defines one-parameter semigroups $\{p^{c,d}\}_{c,d\in\kil}$
in the unital Banach algebra $(\Bil^*,\conv)$ which we refer to as the
\emph{associated convolution semigroups} of the cocycle. The cocycle is
said to be \emph{Markov-regular} if each of its associated semigroups is
norm-continuous.

\begin{defn}
A process $l\in \Procstar(\Bil)$ is a
(\emph{weak}) \emph{QS convolution cocycle on} $\Bil$
if its associated family
$\{
\omega_{\varepsilon(f_{[0,t[}), \varepsilon(f'_{[0,t[})} \circ l_t |
\ f,f'\in\Step, t\in\Rplus
\}$
is a form QS cocycle on $\Bil$.
\end{defn}

\begin{rems}
(i)
Let $l\in\Procstar(\Bil)$ be a completely bounded QS convolution cocycle on
$\Bil$.
Then $l$ is a QS convolution cocycle in the full sense:
\[
l_{s+t} = l_s \conv \big( \sigma_s \circ l_t \big), \quad
l_0 = \iota_\FFock \circ \counit, \quad
s,t\in\Rplus,
\]
where $(\sigma_s)_{s\geq 0}$ is the injective *-homomorphic
semigroup of right shifts on $B(\FFock)$ and the identification
\[
B(\FFock) =
B(\FFock_{[0,s[})
\otol
\sigma_s\big( B(\FFock)\big)
\]
is invoked.

(ii)
It follows from the proof of Theorem~\ref{Theorem: SLQSDE} that,
for $\varphi \in SL\big(\khat,\khat; \Bil^* \big)$,
the unique form solution of the QS differential equation
\begin{equation}
\label{Equation: bialg QSDE}
dl_t = l_t \conv d\Lambda_\varphi (t), \quad
l_0 = \iota_{\FFock} \circ \counit ,
\end{equation}
is a Markov-regular weak QS convolution cocycle on $\Bil$.

(iii)
Form-cocycles may equally be defined on $\Biltilde$ and
$\Bilol$ with the requirement of strictness/normality,
and $\counit$ replaced by $\counittilde$, respectively $\counitol$.
From the correspondence~\eqref{diag CB} it follows that any one of these
uniquely determines the others.
\end{rems}

Our essential strategy for analysing QS convolution cocycles is to
work in the universal enveloping von Neumann bialgebra $\Bilol$ and,
by transferring between convolution and standard QS cocycles using
the maps $\Rbar$ and $\Ebar$, to apply the theory developed in
[$\text{LW}_{\!1\text{-}4}$], and~\cite{LSqsde}.

We first establish a converse to Remark (ii) above.

\begin{propn}
\label{Proposition: YY}
Let $l$ be a Markov-regular, completely positive, contractive quantum
stochastic convolution cocycle on $\Bil$.
Then there is a unique map $\varphi \in CB(\Bil; B(\khat))$ such
that $l = l^\varphi$.
\end{propn}
\begin{proof}
Set $k := \big(
\Rbar
\ol{l}_t \big)_{t\geq 0}$,
where $\ol{l}_t := \ol{l_t}$ ($t\geq 0$).
Then $k$ is a standard quantum stochastic cocycle on $\Bilol$ which is
Markov-regular, completely positive, contractive and normal.
Therefore, by Theorem 5.10 of~\cite{LWjfa} and Theorem 5.3 of~\cite{LWptrf},
$k$ has a stochastic generator
$\ol{\phi} \in \CBuw\big(\Bilol; \Bilol \otol B(\khat) \big)$,
moreover for $c, d \in \kil$, its associated semigroup $P^{c,d}$ has
generator
$(\id_{\Bilol}\otol \omega_{\chat,\dhat}) \circ \ol{\phi}$.
Set
$\varphiol :=
\Ebar
\ol{\phi} \in
\CBuw\big(\Bilol; B(\khat) \big)$.
Since $\ol{l}_t =
\Ebar
k_t$, the associated
convolution semigroup $p^{c,d}$ of $\ol{l}$ has generating functional
\[
\counitol \circ
\big( \id_{\Bilol} \otol \omega_{\chat,\dhat}\big) \circ \ol{\phi}
=
\omega_{\chat,\dhat} \circ \ol{\varphi}
\]
which equals the generating functional of
the associated convolution semigroup of the QS convolution cocycle
$\ol{l}^\varphi$.
It follows that $\ol{l} = \ol{l}^\varphi$
where $\varphi:= \varphiol|_\Bil$ and so $l = l^\varphi$.
\end{proof}

We refer to $\varphi$ as the \emph{stochastic generator} of the
QS convolution cocycle $l$.
The proof of the next result now proceeds similarly
to those of Theorems 5.1 and 6.1 in~\cite{QSCC2}.

\begin{thm}
\label{Theorem: cocycle}
Let $l$ be  Markov-regular quantum stochastic convolution cocycle on $\Bil$.
Then the following equivalences hold:
\begin{alist}
\item
\begin{rlist}
\item
$l$ is completely positive and contractive\tu{;}
\item
$l = l^\varphi$ for a map $\varphi\in CB(\Bil; B(\khat))$ which is
expressible in the form $\varphi_1 - \varphi_2$ where
$\varphi_1 \in CP(\Bil; B(\khat))$ and $\varphi_2 = \counit(\cdot)
\big( \dQS + |\zeta\ra\la e_0| + |e_0\ra\la \zeta| \big)$
for some $\zeta\in\khat$, and satisfies $\varphitilde (1) \leq 0$.
\end{rlist}
\noindent In this case, the convolution cocycle $l$ is preunital if
and only if its stochastic generator $\varphi$ satisfies
$\varphitilde(1) = 0$.
\item
\begin{rlist}
\item
$l$ is completely positive and preunital\tu{;}
\item
$l = l^\varphi$ for a map $\varphi \in CB(\Bil; B(\khat))$
expressible in the form
\begin{equation}
\label{Equation: instead+} a \mapsto
\begin{bmatrix} \la \xi | \\ D^* \end{bmatrix}
\nu( a )
\begin{bmatrix} | \xi \ra & D \end{bmatrix}
\text{ for } \nu := \rho - \iota_\Kil \circ \counit,
\end{equation}
in which $(\rho, \Kil)$ is a nondegenerate *-representation of
$\Bil$ \tu{(}as $C^*$-algebra\tu{)}, $D$ is an isometry in
$B(\kil;\Kil)$ and $\xi$ is  a vector in $\Kil$.
\end{rlist}
\item
\begin{rlist}
\item
$l$ is *-homomorphic\tu{;}
\item
$l = l^\theta$ where $\theta$ is an $\counit$-structure map\tu{;}
\item
$l = l^\theta$ for a map $\theta$ expressible in the form
\begin{equation}
\label{Equation: instead++} a \mapsto
\begin{bmatrix} \la c | \\ I_\kil \end{bmatrix}
\nu(a)
\begin{bmatrix} | c \ra & I_\kil \end{bmatrix}
\text{ where } \nu := \pi - \iota_\kil \circ \counit,
\end{equation}
for a *-homomorphism $\pi: \Bil \to B(\kil)$ and vector $c\in\kil$.
\end{rlist}
\noindent
In this case, the convolution cocycle $l$ is nondegenerate
if and only if the *-representation $\pi$ is.
\end{alist}
\end{thm}

\begin{proof}
In case (i) of (a), (b) and (c) we let $\varphi$ be the stochastic
generator of $l$, let
$\ol{l} = \ol{l}^\varphi = \big( \ol{l_t} \big)_{t\geq 0}$, and set
$k = k^{\ol{\phi}}$ where
$\ol{\phi} =
\Rbar
\varphiol \in
\CBuw\big( \Bilol; \Bilol \otol B(\khat) \big)$.
Thus $k$ is a Markov-regular standard QS cocycle on $\Bilol$ and
$\varphiol =
\Ebar
\ol{\phi}$.

(a)
If (i) holds then $k$ is completely positive and contractive,
by Proposition~\ref{Proposition: XX},
and normal. Therefore, by Theorem 5.10 of~\cite{LWjfa}, there is a map
$\Phi \in \CPuw( \Bilol; \Bilol \otol B(\khat))$ and operator
$Z \in \Bilol \otol \la \khat |$ such that
\begin{equation}
\label{Equation: phibar and Z}
\ol{\phi}(x) = \Phi(x) -
\big( x \ot \dQS + Z^* (x \ot \la e_0 |) + (x \ot | e_0 \ra) Z \big)
\quad
(x \in \Bilol)
\end{equation}
and $\ol{\phi}(1)\leq 0$.
It follows that $\varphiol(1)\leq 0$ and
\[
\varphiol =
\Psi -
\counitol(\cdot)
\big( \dQS + | \zeta \ra \la e_0| + | e_0 \ra \la \zeta | \big)
\]
where
$\Psi =
\Ebar
\Phi$
and
$\la \zeta | =
\big( \counitol \otol \id_{\la \khat |} \big)(Z)$.
Thus (ii) holds with $\psi= \Psi|_{\Bil}$,
moreover if $l$ is preunital then $\ol{l}$ is unital and so $k$ is
too, therefore $\ol{\phi}(1)=0$ so $\varphiol(1)=0$ also.
Conversely, if (ii) holds then, taking normal extensions,
\[
\varphiol =
\ol{\psi} -
\counitol(\cdot)
\big( \dQS + | \zeta \ra \la e_0| + | e_0 \ra \la \zeta | \big)
\]
and so~\eqref{Equation: phibar and Z} holds with
$\Phi =
\Rbar
\ol{\psi}$ and
$Z = 1_{\Bilol} \ot \la \zeta |$. Therefore,
by~\cite{LWptrf} Theorem 5.3,
$k$ is completely positive and contractive and so, by
Proposition~\ref{Proposition: XX}, $l$ is too. Similarly, if
$\varphiol(1)=0$
then
$\ol{\phi}(1)=0$ so $k$ is unital,
thus $\ol{l}$ is too,
and therefore $l$ is preunital. This proves (a).

(b)
If (i) holds then, choosing $\psi$ and $\zeta$ as in (a), let
\[
\begin{bmatrix} \la \xi | \\ D^* \end{bmatrix}
\rho( \cdot )
\begin{bmatrix} | \xi \ra & D \end{bmatrix}
\]
be a minimal Stinespring decomposition of $\psi$.
Thus $(\rho, \Kil)$ is a nondegenerate representation of $\Bil$ and
\[
\big( \dQS + |\zeta\ra\la e_0| + |e_0\ra\la \zeta| \big)
= \wt{\psi}(1) =
\begin{bmatrix}
\|\xi\|^2 & \la \xi | D \\
D^* | \xi \ra & D^*D
\end{bmatrix},
\]
so $D$ is isometric and (ii) holds. Conversely, suppose that (ii)
holds. Then $\varphiol(1)=0$ and it is easily verified that
$\varphi$ has the form given in Part (ii) of (a), with
$\zeta = \binom{\frac{1}{2}\|\xi\|^2}{D^*\xi}$, therefore (i) holds
by (a). This proves (b).

(c)
If (i) holds then $k$ is *-homomorphic so,
by~\cite{LWptrf} Proposition 6.3,
$\ol{\phi}$ is a structure map:
\begin{equation}
\label{Equation: (b)} \ol{\phi} (x^*y) = \ol{\phi}(x)^* \iota(y) +
\iota(x)^* \ol{\phi}(y) + \ol{\phi} (x)^*  (1_{\Bilol} \ot \dQS)
\ol{\phi} (y) \quad (x,y\in \Bilol),
\end{equation}
where $\iota$ denotes the ampliation map $x \mapsto x \ot I_\FFock$.
Since
$\counitol \otol \id_{B(\khat)}$ is a unital *-homomorphism this implies
that
\begin{equation}
\label{Equation: varphibar structure}
\varphiol (x^*y) =
\varphiol(x)^* \counitol (y) + \counitol(x)^* \varphiol(y) +
\varphiol (x)^* \dQS \,\varphiol (y)
\quad
(x,y\in\Bilol)
\end{equation}
and so (ii) holds.
Suppose conversely that (ii) holds. By separate $\sigma$-weak continuity
of multiplication in $\Bilol$ it follows
that~\eqref{Equation: varphibar structure} holds and a brief calculation
confirms the identity
\[
\Omega(u^*v) =
\big( \id_{\Bilol} \otol \counitol \big)(u)^*
\Omega(v) +
\Omega(u)^*
\big( \id_{\Bilol} \otol \counitol \big)(v) +
\Omega(u)^*
\big( 1_{\Bilol} \ot \dQS \big)
\Omega(v),
\]
where $\Omega := \big( \id_{\Bilol} \otol \varphiol \big)$,
for simple tensors $u,v\in\Bilol\otul\Bilol$. Since both sides are
separately $\sigma$-weakly continuous the identity is valid for all $u$
and $v$ in $\Bilol \otol \Bilol$.
Substituting in $u = \coproductol x$ and $v = \coproductol y$ we see that
$\ol{\phi}$ satisfies~\eqref{Equation: (b)}.
Therefore, by Corollary 4.2 of
by~\cite{LWhomomorphic},
$k$ is *-homomorphic thus, by
Proposition~\ref{Proposition: XX}, $l$ is too
and therefore (ii) holds. The equivalence of
(ii) and (iii) is the general form of an $\counit$-structure map
(see~\eqref{Equation: block}).
In view of (a), the last part is easily seen from the representation (iii).
This completes the proof.
\end{proof}

\begin{rem} 
The proper hypothesis for Parts (a) and (b) above is that $\Bil$ be a
(multiplier) $\Cstar$-\emph{hyperbialgebra}, since the multiplicative
property of $\coproduct$ is not used in their proof.
The above result therefore generalises Theorems 5.1 and 6.2
of~\cite{QSCC2} to the locally compact category.
\end{rem}

\section{Quantum L\'evy processes on multiplier $C^*$-bialgebras}
\label{Section: Levy}

In this section we extend the definition of weak quantum L\'evy process to
multiplier $C^*$-bialgebras and establish a reconstruction theorem which is
analogous to Sch\"urmann's for purely algebraic bialgebras (\cite{Schurmann})
and extends ours, proved for unital $C^*$-bialgebras in \cite{QSCC2}.


Throughout this section $\Bil$ \emph{denotes a fixed counital $C^*$-bialgebra}.

\begin{defn}  \label{wLp}
A \emph{weak quantum L\'{e}vy process} on $\Bil$ over a $C^*$-algebra-with-a-state $(\Al , \omega)$
is a family $\big(j_{s,t}\! :\Bil \to \Altilde\big)_{0 \leq s \leq t}$ of nondegenerate
*-homomorphisms for which the functionals  $\lambda_{s,t}:= \omega \circ j_{s,t}$
satisfy the following conditions, for $0\leq r \leq s \leq t$:

\begin{rlist}
\item
$\lambda_{r,t} = \lambda_{r,s} \star \lambda_{s,t}$\tu{;}
\item
$\lambda_{t,t} = \Cou $\tu{;}
\item
$\lambda_{s,t} = \lambda_{0,t-s}$\tu{;}
\item
\[
\wt{\omega} \left( \prod^n_{i=1} j_{s_i,t_i} (x_i) \right) = \prod^n_{i=1}
\lambda_{s_i,t_i} (x_i)
\]
whenever $n \in \bn$, $x_1, \ldots, x_n \in\Bil$ and
the intervals $[s_1,t_1[,\ldots ,[s_n, t_n[$ are disjoint\tu{;}
\item
$\lambda_{0,t}  \to \Cou$ pointwise as $t \to 0$.
\end{rlist}
A weak quantum L\'{e}vy process is called
\emph{Markov-regular} if
$\lambda_{0,t} \to\Cou$ in norm, as $t \to 0$.
\end{defn}

\begin{rems}
In the case of unital $C^*$-bialgebras we did not insist
that the *-algebra $\Al$ was a $C^*$-algebra.

As in the unital case, we refer to the weakly continuous convolution semigroup
$(\lambda_t := \lambda_{0,t})_{t \geq 0}$ on $\Bil$ as the
\emph{one-dimensional distribution} of the process, and call the process
\emph{Markov-regular} if this is norm-continuous, in which case we refer to the
convolution semigroup generator as the \emph{generating functional} of the process (\cite{discrete}).
Moreover, as in the unital case, we call two weak quantum L\'evy processes \emph{equivalent}
if their one-dimensional distributions coincide.
\end{rems}


The generating functional $\gamma$ of a Markov-regular weak quantum L\'evy process,
being the generator of a norm-continuous convolution semigroup of states,
is \emph{real},
that is $\gamma = \gamma^\dagger$ where $\gamma^\dagger(a):= \overline{\gamma(a^*)}$,
\emph{conditionally positive},
that is positive on the ideal $\Ker \counit$,
and its strict extension satisfies
$\gammatilde (1) = 0$.
Note that if $l\in \Procstar(\Bil)$ is a QS convolution cocycle on $\Bil$,
with noise dimension space $\kil$,
which is *-homomorphic and preunital then, setting
$\Al := K(\FFock)$, $\omega:= \omega_{\ve(0)}$,
and $j_{s,t} := \sigma_s \circ l_{t-s}$ for all $0\leq s\leq t$,
we obtain a weak quantum L\'evy process on $\Bil$,
called a \emph{Fock space quantum L\'evy process}, which is Markov-regular if $l$ is.
Our goal now is to establish a converse, in other words to extend the
reconstruction theorem of~\cite{QSCC2} to the nonunital case.
We give an elementary self-contained proof, independent of
automatic implementability/complete boundedness properties of
$\chi$-structure maps.
Recall Lemma~\ref{condpos}.

\begin{thm} \label{Crecon}
Let $\gamma\in\Bil^*$ be real,
conditionally positive and satisfy $\wt{\gamma}(1)=0$.
Then there is a \tu{(}Markov-regular\tu{)}
Fock space quantum L\'{e}vy process with generating functional $\gamma$.
\end{thm}

\begin{proof}
By Theorem~\ref{Theorem: cocycle} it suffices to show that there is a
Hilbert space $\kil$ and an
$\counit$-structure map $\varphi: \Bil \to B(\kilhat)$ of the form
$\left[\begin{smallmatrix}\gamma & *\\ * & * \end{smallmatrix}\right]$
satisfying $\varphitilde (1) = 0$.
Set $\gammatilde_0 := \gammatilde|_{\Ker \counittilde}$ and let
$\psi$ be the map
$\Bil \to \Biltilde$, $b \mapsto b - \counit (b)1$.
By Theorem~\ref{Theorem 1.1} and Lemma~\ref{condpos},
$\gammatilde$ is real and $\gammatilde_0$ is positive.
Since also $\gammatilde (1) = 0$,
\begin{equation} \label{7.2i}
q:(a,b) \mapsto
\gamma(a^*b) - \gamma(a)^*\counit(b) - \counit(a)^*\gamma(b) =
\gammatilde_0 \big(\psi(a)^*\psi(b)\big)
\end{equation}
defines a nonnegative sesquilinear form on $\Bil$.
Let $\kil$ and $d: \Bil \to \kil$ be respectively
the Hilbert space and induced map obtained by
quotienting $\Bil$ by the null space of $q$ and completing,
so that
\[
\ol{d(\Bil)} = \kil \text{ and } \la d(a), d(b) \ra = q(a,b), \quad a,b \in \Bil,
\]
and let $\delta$ be the linear map
$\Bil \to | \kil \ra$, $b\mapsto |d(b)\ra$.
Then, by the complete boundedness of $\gammatilde$ and $\psi$,
\[
\big\| \delta^{(n)}(A) u \big\|^2 =
\big\la u,
(\gammatilde_0)^{(n)}\big(\psi^{(n)}(A)^*\psi^{(n)}(A) \big) u \big\ra
\leq
\|\gammatilde\|_{\cb}\|\psi\|^2_{\cb} \|A\|^2\|u\|^2,
\]
for all $n\in\mathbb{N}$, $A\in M_n(\Bil)$ and $u \in \Comp^n$,
so $\delta$ is completely bounded and we have
\begin{equation}
\label{7.A}
\delta(a)^*\delta(b) =
\gamma(a^*b) - \gamma(a)^*\counit(b) - \counit(a)^*\gamma(b),
\quad a,b\in\Bil.
\end{equation}
Now
\begin{align*}
\big\| d(ab) - \counit(b)d(a) \big\|^2 &=
\gammatilde_0 \big( \psi(b)^* a^*a \psi(b) \big)
\\ &\leq
\| a \|^2 \,\gammatilde_0\big( \psi(b)^*\psi(b)\big) = \| a \|^2 \| d(b) \|^2,
\quad a,b \in \Bil,
\end{align*}
so there are bounded operators $\pi (a)$ on $\kil$ satisfying
\begin{equation} \label{7.2ii}
\pi(a)d(b) = d(ab) - \counit(b)d(a), \quad a,b\in\Bil.
\end{equation}
Using the density of $d(\Bil)$ it is straightforward to verify that
the map $a\mapsto \pi(a)$ defines a *-representation of $\Bil$ on $\kil$.
From~\eqref{7.2ii},
$\delta$ is a $(\pi,\counit)$-derivation and so,
from~\eqref{7.A},
$\varphi :=
\left[ \begin{smallmatrix}
\gamma & \delta^\dagger \\ \delta & \pi - \iota_{\kil}
\end{smallmatrix} \right]$
defines an $\counit$-structure map $\Bil \to B(\kilhat)$,
and therefore it only remains to prove that $\varphitilde(1)=0$.
Since $\gammatilde(1)=0$, this follows from the identities
\[
\delta(a)^*\delta(b) =
\gamma\big(a^*b - \counit(a)^*b - a^*\counit(b)\big)
\text{ and }
\pi(a)\delta(b) = \delta(ab) - \delta(a)\counit(b),
\quad a,b\in\Bil,
\]
and the density of $\bigcup\{\Ran \delta(b): b\in\Bil\} = d(\Bil)$ in $\kil$.
\end{proof}


This has two significant consequences.

\begin{cor}
Every Markov-regular weak quantum L\'{e}vy process is equivalent to a Fock
space quantum L\'{e}vy process.
\end{cor}

The second consequence uses the deeper fact that every
$\counit$-structure map is implemented (see Theorem~\ref{established}).

\begin{thm}
\label{significant}
Let $\gamma \in \Bil^*$. Then the following are equivalent:
\begin{rlist}
\item
$\gamma$ is the generating functional of a norm-continuous
convolution semigroup of states on $\Bil$\tu{;}
\item
$\gamma$ is real, conditionally positive and
satisfies $\gammatilde(1)=0$\tu{;}
\item
There is a nondegenerate representation $(\pi,\hil)$ of $\Bil$
and vector $\eta\in\hil$ such that
$\gamma = \omega_\eta \circ (\pi - \iota_\hil \circ \counit)$.
\end{rlist}
\end{thm}

As stated earlier, the above results mean that all the Poisson-type
convolution semigroups of states on $\Bil$ are easily constructed,
along with their associated quantum L\'evy processes.

In~\cite{QSCC2} we also introduced a stronger notion of
\emph{product system quantum L\'evy processes} on a unital and
counital $C^*$-bialgebra $\Bil$ and established the following two
facts: each Fock space quantum L\'evy process on $\Bil$ is in
particular a product system quantum L\'evy process and each product
system quantum L\'evy process determines in a natural way a weak
quantum L\'evy process on $\Bil$ with the same finite-dimensional
distribution. The definition of a product system quantum L\'evy
process extends naturally to the nonunital case, with the assumption
of \emph{unitality} --- of the $^*$-homomorphisms constituting the
process --- replaced by \emph{nondegeneracy}, and the proofs of the
above two facts remain valid.



\end{document}